\documentclass{amsart}
\usepackage{pstricks,amssymb}

\title{Trees, parking functions, syzygies, 
and deformations of monomial ideals}
\author{Alexander Postnikov \and Boris Shapiro}

\address{Department of Mathematics, 
M.I.T., Cambridge, MA 02139, U.S.A.}
\email{apost@math.mit.edu}
\urladdr{http://www.math-mit.edu/\mytilde apost/}

\address{Department of Mathematics, 
University of Stockholm, Stockholm, S-10691, Sweden}
\email{shapiro@matematik.su.se} 
\urladdr{http://www.matematik.su.se/\mytilde shapiro/}

\date{December~17, 2002; revised on February~2, 2003;
revised on February~4, 2003}


\thanks{The first author was supported in part by NSF grant DMS-0201494.}


\newtheorem{theorem}{Theorem}[section]
\newtheorem{proposition}[theorem]{Proposition}
\newtheorem{corollary}[theorem]{Corollary}
\newtheorem{lemma}[theorem]{Lemma}
\newtheorem{conjecture}[theorem]{Conjecture}
\newtheorem{example}[theorem]{Example}

\def\<{\left<}
\def\>{\right>}
\def\({\left(}
\def\){\right)}
\def\A{\mathcal{A}}
\def\B{\mathcal{B}}
\def\C{\mathcal{C}}
\def\S{\mathcal{S}}
\def\M{\mathcal{M}}
\def\K{\mathbb{K}}
\def\Z{\mathbb{Z}}
\def\I{\mathcal{I}}
\def\J{\mathcal{J}}
\def\R{\mathcal{R}}
\def\MF{\mathfrak{M}}
\def\MM{M^{\textrm{\tiny max}}}
\def\Ker{\mathrm{Ker}}
\def\Hilb{\mathrm{Hilb}}
\def\mytilde{\kern-.015in\hbox{\lower.03in\hbox{\~{}}}\kern-.01in}
\def\ds{\displaystyle}
\def\lcm{\mathrm{lcm}}

\psset{unit=.5pt, linewidth=.5pt}
\newrgbcolor{mygreen}{.2 .7 .2}

\begin{document}

\begin{abstract}
For a graph $G$, we construct two algebras, whose dimensions 
are both equal to the number of spanning trees of $G$.
One of these algebras is the quotient of the polynomial ring
modulo certain monomial ideal, while the other is the quotient 
of the polynomial ring modulo certain powers of linear forms.
We describe the set of monomials that forms a linear basis in each 
of these two algebras.  The basis elements correspond to 
$G$-parking functions that naturally came up in the abelian
sandpile model.
These ideals are instances of the general class of monotone monomial
ideals and their deformations.  We show that the Hilbert series of a
monotone monomial ideal is always bounded by the Hilbert series of its
deformation.  Then we define an even more general class of monomial ideals 
associated with posets and construct free resolutions for these ideals.
In some cases these resolutions coincide with Scarf resolutions.
We prove several formulas for Hilbert series of monotone monomial ideals 
and investigate when they are equal to Hilbert series of deformations.
In the appendix we discuss the abelian sandpile model.
\end{abstract}

\maketitle
\pagestyle{myheadings}
\markboth{\sc alexander postnikov and boris shapiro}
{\sc trees, parking functions, syzygies, and deformations of monomial ideals}

\section{Introduction}

The famous formula of Cayley says that the number of trees on $n+1$ labelled
vertices equals $(n+1)^{n-1}$.  Remarkably, this number has several other
interesting combinatorial interpretations.  For example, it is equal to the
number of parking functions of size $n$.

In this paper we present two algebras $\A_n$ and $\B_n$ of dimension
$(n+1)^{n-1}$.  The algebra $\A_n$ is a quotient of the polynomial ring
modulo a monomial ideal; and the algebra $\B_n$ is a quotient of the
polynomial ring modulo some powers of linear forms.  It is immediate that
the set of monomials $x^b$, where $b$ is a parking function, is the standard
monomial basis of the algebra $\A_n$.  On the other hand, the same set of
monomials forms a basis of the algebra $\B_n$, which is a non-trivial
result.

More generally, for any graph $G$, we define two algebras $\A_G$ and $\B_G$
and describe their monomial bases.  
The basis elements correspond to $G$-parking functions.
These functions extend the usual parking functions and are related
to the abelian sandpile model;
their number  equals the number of spanning trees of the graph $G$.  
This implies that $\dim\A_G=\dim\B_G$ is also the 
number of spanning trees of $G$.

All these pairs of algebras are instances of the general class of algebras
given by monotone monomial ideals and their deformations.  For such an
algebra $\A$ and its deformation $\B$, we show that $\dim \A\geq \dim\B$ and
the Hilbert series of $\B$ is termwise bounded by the Hilbert series of
$\A$.  There is a natural correspondence between polynomial generators of
the ideal for $\B$ and monomial generators of the ideal for $\A$.  However,
these monomials are not the leading terms of the polynomial generators for
any term order, because they are usually located at the center of the Newton
polytope of the corresponding polynomial generators.    The standard
Gr\"obner bases technique cannot be applied to this class of algebras.

We also investigate the class of order monomial ideals that extends monotone
monomial ideals.  These ideals are associated with posets whose elements are
marked by monomials.  We construct a free resolution for such an ideal as the
cellular resolution corresponding to the order complex of the poset.  This
resolution is minimal if the ideal satisfies some generosity condition.  In
this case, the numbers of increasing $k$-chains in the poset are exactly the
Betti numbers of the ideal.  This resolution often coincides with the Scarf
resolution.

We discuss some results of our previous works on the algebra generated by the
curvature forms on the generalized flag manifold.  This algebra extends
the cohomology ring of the generalized flag manifold.  
For type $A_{n-1}$, the dimension this algebra
equals the number of forests on $n+1$ vertices.  The algebras generated by the
curvature forms are analogous to the algebras that we study in the present
paper.  This attempt to lift Schubert calculus on the level of differential
forms was our original motivation.

\medskip
The general outline of the paper follows.  In Section~\ref{sec:G-parking} we
define $G$-parking functions for a digraph $G$.  We formulate
Theorem~\ref{th:G-parking=trees} that says that the number of such functions
equals the number of oriented spanning trees of $G$.  Then we construct the
algebra $\A_G$ as the quotient of the polynomial ring modulo certain monomial
ideal.  Elements of the standard monomial basis of $\A_G$ correspond
to $G$-parking functions.
In Section~\ref{sec:power} we construct the
algebra $\B_G$ as the quotient of the polynomial ring modulo the ideal
generated by power of certain linear forms.  Then we formulate
Theorem~\ref{th:ABK} that implies that the algebras $\A_G$ and $\B_G$ have
the same Hilbert series.  In Section~\ref{sec:examples} we give two examples
of these results.  For the complete graph $G=K_{n+1}$ we recover the usual
parking functions and the algebras $\A_n$ and $\B_n$ of dimension
$(n+1)^{n-1}$.  For a slightly more general class of graphs we obtain two
algebras of dimension $l\,(l+kn)^{n-1}$.  Section~\ref{sec:monotone} is
devoted to description of monotone monomial ideals and their deformations.
We formulate Theorem~\ref{th:span}, which implies the inequality for the
Hilbert series.  In Section~\ref{sec:min-free-res} we describe a more
general class of monomial ideals associated with posets and construct free 
resolutions for these ideals.  Components of the resolution for such an 
ideal correspond to strictly increasing chains in the poset.
In Section~\ref{sec:three-examples} we give several examples of minimal
free resolutions.
In Section~\ref{sec:hilbert-monomial} we prove general
formulas for the Hilbert series and dimension of the algebra given by a monotone
monomial ideal.  Then we deduce Theorem~\ref{th:G-parking=trees}.  
In Section~\ref{sec:CG} we construct the
algebra $\C_G$ and prove Theorem~\ref{th:CG} that claims that the dimension
of this algebra equals the number of spanning trees.  Actually, we will
later see that $\C_G$ is isomorphic to the algebra $\B_G$.  In
Section~\ref{sec:monote-proof} we prove Theorem~\ref{th:span}.  Then we
finish the proof of Theorem~\ref{th:ABK}, which goes as follows.  By
Theorem~\ref{th:span} and construction of $\C_G$ we know that
$\Hilb\,\A_G\geq \Hilb\,\B_G\geq \Hilb\,C_G$ termwise.  On the other hand, by
Theorems~\ref{th:G-parking=trees} and~\ref{th:CG}, $\dim\A_G=\dim\C_G$ is
the number of spanning trees of $G$.   Thus the Hilbert series of these
three algebras coincide.  In Section~\ref{sec:forests} we discuss some
results of our previous works and compare them with results of this paper.
We mention a certain algebra, whose dimension equals the number of forests
on $n+1$ vertices.  This algebra originally appeared in the attempt to lift
Schubert calculus of the flag manifold on the level of differential forms.
In Section~\ref{sec:rho-algebras} we discuss a special class of monotone
monomial ideals and their deformations.  We give a minimal free resolution
and subtraction-free formula for the Hilbert series of the algebra $\A$ 
and list several cases when it is equal to the Hilbert series of $\B$.
The appendix is devoted to the abelian sandpile model and its links
with $G$-parking functions.

\bigskip
{\sc Acknowledgments:}
We are  grateful to Richard Stanley,  Mikhail Shapiro, Bernd
Sturmfels, Ezra Miller, Hal Schenck, Gilles Schaeffer, 
Andrei Gabrielov,
and Andrei Zelevinsky for helpful discussions and relevant
comments, and to Ralf Fr\"oberg for his help with {\tt Macaulay2}.  The first
author was supported by the Miller Institute 
%for Basic Research in Science 
at UC Berkeley during 1999--2001,  when a part of this project was
completed, and by NSF grant DMS-0201494.
The second author is sincerely grateful to the Max-Planck Institut
f\"ur Mathematik in Bonn for the financial support and the excellent research
atmosphere during his visit in 2000.

\section{$G$-parking functions}
\label{sec:G-parking}

A {\it parking function} of size $n$ is a sequence $b=(b_1,\dots,b_n)$ of
non-negative integers such that its increasing rearrangement $c_1\leq
\cdots\leq c_n$ satisfies $c_i< i$.  Equivalently, we can formulate this
condition as $\#\{i\mid b_i < r\} \geq r$, for $r=1,\dots,n$.  The parking
functions of size $n$ are known to be in bijective correspondence with
trees on $n+1$ labelled vertices, see Kreweras~\cite{Krew}.  
Thus, according to Cayley's formula for the number of labelled trees,
the total number of parking functions of size $n$ equals $(n+1)^{n-1}$. 
In this section we extend this statement to a more general class of functions.  
\medskip

A {\it graph\/} is given by specifying its set of vertices, set of edges,
and a function that associates to each edge an unordered pair of vertices.
A {\it directed graph}, or {\it digraph}, is given by specifying its set of
vertices, set of edges, and a function that associates to each edge an
ordered pair of vertices.  Thus multiple edges and loops are allowed in
graphs and digraphs.  A {\it subgraph\/} $H$ in a (directed) graph $G$ is a
(directed) graph on the same set of vertices whose set of edges is a subset
of edges of $G$.  We will write $H\subset G$ to denote that $H$ is a subgraph
of $G$.  For a subgraph $H\subset G$, let $G\setminus H$ denote the
{\it complement subgraph\/} whose edge set is complementary to that of $H$.
Also we will write $e\in G$ to show that $e$ is an edge of the graph $G$.

Let $G$ be a digraph on the set of vertices $0,1,\dots,n$.  The vertex $0$
will be the root of $G$.  The digraph $G$ is determined by its
{\it adjacency matrix} $A=(a_{ij})_{0\leq i,j\leq n}$, where $a_{ij}$ is the
number of edges from the vertex $i$ to the vertex $j$.  We will regard
graphs as a special case of digraphs with symmetric adjacency matrix $A$.

An {\it oriented spanning tree\/} $T$ of the digraph $G$ is a subgraph 
$T\subset G$ such that there exists a unique directed path in $T$ 
from any vertex $i$ to the root $0$.  The number $N_G$ of such trees
is given by the {\it Matrix-Tree Theorem}, e.g., see~\cite[Section~5.6]{EC2}:
\begin{equation}
N_G = \det L_G,
\label{eq:matrix-tree}
\end{equation}
where $L_G=(l_{ij})_{1\leq i,j\leq n}$ the 
{\it truncated Laplace matrix},
also known as the {\it Kirkhoff matrix}, given by
\begin{equation}
l_{ij} = \left\{
\begin{array}{cl}
\ds \sum_{r\in\{0,\dots,n\}\setminus\{i\}} a_{ir}
& \textrm{for } i = j, \\[.1in]
\ds -a_{ij}            & \textrm{for } i\ne j.
\end{array}
\right.
\label{eq:Laplace}
\end{equation}
If $G$ is a graph, i.e., $A$ is a symmetric matrix, then oriented spanning
trees defined above are exactly the usual {\it spanning trees\/} of $G$,
which are connected subgraphs of $G$ without cycles.

For a subset $I$ in $\{1,\dots,n\}$ and a vertex $i\in I$, let 
$$
d_{I}(i) = \sum_{j\not\in I} a_{ij},
$$ 
i.e., $d_I(i)$ is the number of edges from the
vertex $i$ to a vertex outside of the subset $I$.  Let us say that a
sequence $b=(b_1,\dots,b_n)$ of non-negative integers is a {\it $G$-parking
function\/} if, for any nonempty subset $I\subseteq\{1,\dots,n\}$, there
exists $i\in I$ such that $b_i<d_I(i)$.

If $G=K_{n+1}$ is the complete graph on $n+1$ vertices then
$K_{n+1}$-parking functions are the usual parking functions of size $n$
defined in the beginning of this section.

\begin{theorem} {\rm cf.~\cite{Gab1}} \
The number of $G$-parking functions equals the number
$N_G=\det L_G$ of oriented spanning trees of the digraph $G$.
\label{th:G-parking=trees}
\end{theorem}

Interestingly, $G$-parking functions are related to the abelian
sandpile model introduced by Dhar~\cite{Dhar}.  In the appendix we will
discuss the sandpile model and show that
Theorem~\ref{th:G-parking=trees} is essentially equivalent to the
result of Gabrielov~\cite[Eq.~(21)]{Gab1} on sandpiles.  In
Section~\ref{sec:hilbert-monomial} we will prove
Theorem~\ref{th:G-parking=trees} without using the sandpile
model.

We can reformulate the definition of $G$-parking functions in algebraic
terms as follows.  Throughout this paper we fix a field  $\K$.
Let $\I_G = \<m_I\>$ be the monomial ideal in 
the polynomial ring $\K[x_1,\dots,x_n]$ 
generated by the monomials
\begin{equation}
m_I = \prod_{i\in I} x_i^{d_I(i)},
\label{eq:mI}
\end{equation}
where $I$ ranges over all nonempty subsets $I\subseteq\{1,\dots,n\}$.
Define the algebra $\A_G$ as the quotient $\A_G=\K[x_1,\dots,x_n]/\I_G$. 

A non-negative integer sequence $b=(b_1,\dots,b_n)$
is a $G$-parking function if and only if the monomial 
$x^b=x_1^{b_1}\cdots x_n^{b_n}$ is nonvanishing in the algebra $\A_G$.

For a monomial ideal $\I$, the set of all monomials that do not belong to
$\I$ is a basis of the quotient of the polynomial ring modulo $\I$,
called the {\it standard monomial basis.}
Thus the monomials $x^b$, where $b$ ranges over $G$-parking
functions, form the standard monomial basis of the algebra $\A_G$.

\begin{corollary}
The algebra $\A_G$ is finite-dimensional as a linear space over $\K$.
Its dimension is equal to the number of oriented spanning trees 
of the digraph $G$:
$$
\dim \A_G = N_G.
$$
\end{corollary}

For an undirected graph, $G$-parking functions and monomials $m_I$ 
also appeared in a recent paper by Cori, Rossin, and Salvy~\cite{CRS}.

\section{Power algebras}
\label{sec:power}

Let $G$ be an undirected graph on the set of vertices $0,1,\dots,n$.
In this case the dimension of the algebra $\A_G$ is equal to the number
of usual spanning trees of $G$.

For a nonempty subset $I$ in $\{1,\dots,n\}$, let 
$D_I=\sum_{i\in I,\,j\not\in I} a_{ij} = \sum_{i\in I} d_I(i)$ 
be the total number of edges that join some vertex in $I$ with 
a vertex outside of $I$.
For any nonempty subset $I\subseteq\{1,\dots,n\}$, let
\begin{equation}
p_I = \left(\sum_{i\in I} x_i\right)^{D_I}.
\label{eq:pI}
\end{equation}

Let $\J_G=\<p_I\>$ be the ideal in the polynomial ring $\K[x_1,\dots,x_n]$ 
generated by the polynomials $p_I$ for all nonempty subsets $I$.
Define the algebra $\B_G$ as the quotient $\B_G=\K[x_1,\dots,x_n]/\J_G$.

The algebras $\A_G$ and $\B_G$, as well as all other algebras in this paper,
are graded.
For a graded algebra $\A=\A^0\oplus\A^1\oplus\A^2\oplus\cdots$,
the {\it Hilbert series\/} of $\A$ is the formal power series
in $q$ given by
$$
\Hilb\,\A=\sum_{k\geq 0} q^k\,\dim\A^k.
$$

Our first main result is the following statement.

\begin{theorem} 
The monomials $x^b$, where $b$ ranges over $G$-parking
functions, form a linear basis of the algebra $\B_G$.
Thus the Hilbert series of the algebras $\A_G$ 
and $\B_G$ coincide termwise: $\Hilb\,\A_G=\Hilb\,\B_G$.  
In particular, both these algebras
are finite-dimensional as linear spaces over $\K$ and
$$
\dim \A_G = \dim \B_G = N_G
$$
is the number of spanning trees of the graph $G$.
\label{th:ABK}
\end{theorem}

\begin{example}
{\rm 
Let $n=3$ and let $G$ be the graph given by
\begin{equation}
G=
\lower.2in\hbox{\pspicture(-20,-10)(70,67)
\psline{-}(0,0)(50,0)(50,50)(0,50)(0,0)
\psline{-}(0,0)(50,0)(0,50)
\pscircle*(0,0){2}
\pscircle*(50,0){2}
\pscircle*(0,50){2}
\pscircle*(50,50){2}
\rput(60,60){{\footnotesize0}}
\rput(-10,60){{\footnotesize1}}
\rput(-10,-10){{\footnotesize2}}
\rput(60,-10){{\footnotesize3}}
\endpspicture}.
\label{eq:example-graph}
\end{equation}
The graph $G$ has 8 spanning trees:
$$
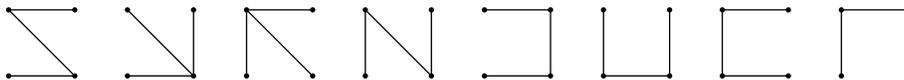
\pspicture(-20,-10)(70,67)
\psline{-}(0,0)(50,0)
\psline{-}(50,50)(0,50)
\psline{-}(50,0)(0,50)
\pscircle*(0,0){2}
\pscircle*(50,0){2}
\pscircle*(0,50){2}
\pscircle*(50,50){2}
\endpspicture
\pspicture(-20,-10)(70,67)
\psline{-}(0,0)(50,0)
\psline{-}(50,50)(50,0)
\psline{-}(50,0)(0,50)
\pscircle*(0,0){2}
\pscircle*(50,0){2}
\pscircle*(0,50){2}
\pscircle*(50,50){2}
\endpspicture
\pspicture(-20,-10)(70,67)
\psline{-}(0,0)(0,50)
\psline{-}(50,50)(0,50)
\psline{-}(50,0)(0,50)
\pscircle*(0,0){2}
\pscircle*(50,0){2}
\pscircle*(0,50){2}
\pscircle*(50,50){2}
\endpspicture
\pspicture(-20,-10)(70,67)
\psline{-}(0,0)(0,50)
\psline{-}(50,50)(50,0)
\psline{-}(50,0)(0,50)
\pscircle*(0,0){2}
\pscircle*(50,0){2}
\pscircle*(0,50){2}
\pscircle*(50,50){2}
\endpspicture
\pspicture(-20,-10)(70,67)
\psline{-}(0,0)(50,0)(50,50)(0,50)
\pscircle*(0,0){2}
\pscircle*(50,0){2}
\pscircle*(0,50){2}
\pscircle*(50,50){2}
\endpspicture
\pspicture(-20,-10)(70,67)
\psline{-}(0,0)(50,0)(50,50)
\psline{-}(0,50)(0,0)
\pscircle*(0,0){2}
\pscircle*(50,0){2}
\pscircle*(0,50){2}
\pscircle*(50,50){2}
\endpspicture
\pspicture(-20,-10)(70,67)
\psline{-}(0,0)(50,0)
\psline{-}(50,50)(0,50)(0,0)
\pscircle*(0,0){2}
\pscircle*(50,0){2}
\pscircle*(0,50){2}
\pscircle*(50,50){2}
\endpspicture
\pspicture(-20,-10)(70,67)
\psline{-}(50,0)(50,50)(0,50)(0,0)
\pscircle*(0,0){2}
\pscircle*(50,0){2}
\pscircle*(0,50){2}
\pscircle*(50,50){2}
\endpspicture
$$
The ideals $\I_G$ and $\J_G$ are given by 
$$
\begin{array}{l}
\I_G=\<x_1^3,x_2^2,x_3^3,x_1^2x_2,x_1^2x_3^2,x_2x_3^2,x_1x_2^0x_3\>,\\[.1in]
\J_G=\<x_1^3,x_2^2,x_3^3,(x_1+x_2)^3,(x_1+x_3)^4, (x_2+x_3)^3, 
(x_1+x_2+x_3)^2\>.
\end{array}
$$
The standard monomial basis of the algebra $\A_G$
is $\{1, x_1, x_2, x_3, x_1^2, x_1x_2, x_2x_3, x_3^2\}$.
The corresponding $G$-parking functions are the exponent
vectors of the basis elements:
$$
(0,0,0),(1,0,0),(0,1,0),(0,0,1),(2,0,0),(1,1,0),(0,1,1),(0,0,2).
$$
We have $\dim\A_G=\dim\B_G=8$ is the number of spanning trees of $G$,
and $\Hilb\,\A_G=\Hilb\,\B_G=1+3q+4q^2$.

}
\end{example}

We will refine Theorem~\ref{th:ABK} and interpret dimensions of graded
components of the algebras $\A_G$ and $\B_G$ in terms of certain statistics
on spanning trees.  Let us fix a linear ordering of all edges of the graph
$G$.  For a spanning tree $T$ of $G$, an edge $e\in G\setminus T$ is called
{\it externally active\/} if there exists a cycle $C$ in the graph $G$ such
that $e$ is the minimal edge of $C$ and $(C\setminus\{e\})\subset T$.  The
{\it external activity\/} of a spanning tree is the number of externally
active edges.  Let $N_G^k$ denote the number of spanning trees $T\subset G$
of external activity $k$.  Even though the notion of external activity
depends on a particular choice of ordering of edges, the numbers $N_G^k$ are
known to be invariant on the choice of ordering.

Let $\A_G^k$ and $\B_G^k$ be the $k$-th graded components of the algebras
$\A_G$ and $\B_G$, correspondingly.

\begin{theorem}  
The dimensions of the $k$-th graded components $\A_G^k$
and $\B_G^k$ are equal to 
$$
\dim \A_G^k = \dim \B_G^k = N_G^{|G|-n-k},
$$
the number of spanning trees of $G$ of external activity
$|G|-n-k$, where $|G|$ denotes the number of edges of $G$.
\label{th:dim-graded}
\end{theorem}

\section{Examples:  tree ideals and their generalizations}
\label{sec:examples}

\subsection{Two algebras of dimension $(n+1)^{n-1}$}
\label{ssec:(n+1)^{n-1}}

Suppose that $G=K_{n+1}$ is the complete graph on $n+1$ vertices.
As we have already mentioned, the $K_{n+1}$-parking functions
are the usual parking functions of size $n$ defined in the beginning
of Section~\ref{sec:G-parking}.

Let $\I_n=\<m_I\>$ and $\J_n=\<p_I\>$ be the ideals in the polynomial 
ring $\K[x_1,\dots,x_n]$ generated by the monomials $m_I$
and the polynomials $p_I$, correspondingly, given by
$$
\begin{array}{c}
\ds m_I=(x_{i_1}\cdots x_{i_r})^{n-r+1},\\[.1in]
\ds p_I=(x_{i_1}+\cdots+x_{i_r})^{r(n-r+1)},
\end{array}
$$
where in both cases $I=\{i_1,\dots,i_r\}$ runs over all
nonempty subsets of $\{1,\dots,n\}$.
Let $\A_n = \K[x_1,\dots,x_n]/\I_n$ and $\B_n = \K[x_1,\dots,x_n]/\J_n$.

\begin{corollary}
\label{cor:A=B}
The graded algebras $\A_n$ and $\B_n$ have the 
same Hilbert series.  They are finite-dimensional,
as linear spaces over $\K$.  Their dimensions are equal to 
$$
\dim\A_n=\dim\B_n=(n+1)^{n-1}. 
$$
The images of the monomials $x^b$, 
where $b$ ranges over parking functions of size $n$, 
form linear bases in both algebras $\A_n$ and $\B_n$.
\end{corollary}

An {\it inversion} in a tree $T$ on the $n+1$ vertices
labelled $0,\dots,n$ is a pair of vertices labelled $i$ and $j$
such that $i>j$ and the vertex $i$ belongs to the shortest path in $T$ 
that joins the vertex $j$ with the root $0$.

\begin{corollary}
The dimension $\dim \A_n^k = \dim \B_n^k$ of the $k$-th graded components 
of the algebras $\A_n$ and $\B_n$ is equal to
\begin{enumerate}
\item[(A)]
the number of parking functions $b$ of size $n$ such that
$b_1+\cdots+b_n = k$;
\item[(B)]
the number of trees on $n+1$ vertices with external activity 
$\binom{n}{2} -k$;
\item[(C)]
the number of trees on $n+1$ vertices with 
$\binom{n}{2} -k$ inversions.
\end{enumerate}
\end{corollary}

It is well known that the numbers (A), (B), and~(C) are equal,
see~\cite{Krew}.  The {\it inversion 
polynomial\/} is defined as the sum
$I_n(q)=\sum_T q^{\# \textrm{ of inversions in }T}$ over
all trees $T$ on $n+1$ labelled vertices.
Thus the Hilbert series of the algebras $\A_n$ and $\B_n$
are equal to
$$
\Hilb\,\A_n=\Hilb\,\B_n = q^{\binom{n}{2}}\,I_n(q^{-1}).
$$

\subsection{Two algebras of dimension $l\,(l+kn)^{n-1}$}
\label{ssec:l(l+kn)^{n-1}}
It is possible to extend the previous example as follows.
Fix two non-negative integers $k$ and $l$.
Let $G=K_{n+1}^{k,l}$ be the complete graph on
the vertices $0,1,\dots,n$ with the edges $(i,j)$, $i,j\ne 0$, 
of multiplicity $k$ and the edges $(0,i)$ of multiplicity $l$.
The $K_{n+1}^{k,l}$-parking functions are the non-negative integer
sequences $b=(b_1,\dots,b_n)$ such that, for $r=1,\dots,n$,
$$
\#\{i\mid b_i< l+k(r-1)\} \geq r.
$$
The definition of these functions
can be also formulated as $c_i<l+(i-1)\,k$, where
$c_1\leq \cdots\leq c_n$ is the increasing rearrangement
of elements of $b$.
Such functions were studied by Pitman and Stanley~\cite{PiSt} and then 
by Yan~\cite{Yan}.
These authors demonstrated that their number
equals $l\,(l+kn)^{n-1}$.  
One can show, using for example the Matrix-Tree 
Theorem~(\ref{eq:matrix-tree}),
that the number of spanning trees in the graph $K_{n+1}^{k,l}$
equals $l\,(l+kn)^{n-1}$.   Thus Theorem~\ref{th:G-parking=trees}
recovers the above formula for the number of $K_{n+1}^{k,l}$-parking 
functions.

Let $\I_{n,k,l}=\<m_I\>$ and $\J_{n,k,l}=\<p_I\>$ 
be the ideals in the ring $\K[x_1,\dots,x_n]$ generated 
by the monomials $m_I$
and the polynomials $p_I$, correspondingly, given by
$$
\begin{array}{c}
\ds m_I=(x_{i_1}\cdots x_{i_r})^{l+k(n-r)},\\[.1in]
\ds p_I=(x_{i_1}+\cdots+x_{i_r})^{r(l+k(n-r))},
\end{array}
$$
where in both cases $I=\{i_1,\dots,i_r\}$ runs over all
nonempty subsets of $\{1,\dots,n\}$.
Let $\A_{n,k,l} = \K[x_1,\dots,x_n]/\I_{n,k,l}$ and 
$\B_{n,k,l} = \K[x_1,\dots,x_n]/\J_{n,k,l}$.

\begin{corollary}
\label{cor:A=Bnkl}
The graded algebras $\A_{n,k,l}$ and $\B_{n,k,l}$ have the 
same Hilbert series.  They are finite-dimensional,
as linear spaces over $\K$.  Their dimensions are
$$
\dim\A_{n,k,l}=\dim\B_{n,k,l}=l\,(l+kn)^{n-1}. 
$$
The images of the monomials $x^b$, 
where $b$ ranges over $K_{n+1}^{k,l}$-parking functions,
form linear bases in both algebras $\A_{n,k,l}$ and $\B_{n,k,l}$.
\end{corollary}

\section{Monotone monomial ideals and their deformations}
\label{sec:monotone}

A {\it monotone monomial family\/} 
is a collection $\M=\{m_I\mid I\in\Sigma\}$ of monomials 
in the polynomial ring $\K[x_1,\dots,x_n]$ labelled by a set $\Sigma$
of nonempty subsets in $\{1,\dots,n\}$ that satisfies the following 
three conditions:

\smallskip
\begin{itemize}
\item[(MM1)]
  For $I\in \Sigma$, $m_I$ is a monomial in the variables $x_i$, $i\in I$.
\smallskip
\item[(MM2)]
  For $I,J\in\Sigma$ such that $I\subset J$ and
$i\in I$, we have 
$\deg_{x_i}(m_I)\geq \deg_{x_i}(m_J)$.
\smallskip
\item[(MM3)]
  For $I,J\in \Sigma$, $\lcm(m_I,m_J)$ is divisible by $m_K$ for
some  $K\supseteq I\cup J$ in $\Sigma$.
\end{itemize}
\smallskip
The {\it monotone monomial ideal\/} $\I=\<\M\>$ associated with 
a  monotone monomial family $\M$ is the ideal in the polynomial ring 
$\K[x_1,\dots,x_n]$ generated by the monomials $m_I$ in $\M$.

It follows from (MM1) and (MM2) that condition (MM3) can be replaced
by the condition:
For $I,J\in \Sigma$ there is $K\supseteq I\cup J$ in $\Sigma$ such that
$m_K$ is a monomial in the $x_i$, $i\in I\cup J$.
This condition is always satisfied if $I,J\in\Sigma$ implies that 
$I\cup J\in\Sigma$.

The monomial ideal $\I_G$ constructed in Section~\ref{sec:G-parking} 
for a digraph $G$ is monotone.  In this case $\Sigma$ is the set of 
all nonempty subsets in $\{1,\dots,n\}$ and $m_I$ is given by~(\ref{eq:mI}).

Remark that two different monotone monomial families may produce
the same monotone monomial ideal.
For example, the ideal $\I_G$, for the graph $G$ shown 
on~(\ref{eq:example-graph}), 
has generator $m_{\{1,3\}}=x_1^2 x_3^2$.
This generator is redundant because it is divisible by 
$m_{\{1,2,3\}}=x_1 x_3$. Thus the same ideal corresponds
to the monotone monomial family with $\Sigma=\{\{1\},
\{2\},\{3\},\{1,2\},\{2,3\},\{1,2,3\}\}$.

Let $I=\{i_1,\dots,i_r\}$.  For a monomial 
$m\in\K[x_{i_1},\dots, x_{i_r}]$, an {\it $I$-deformation} of $m$ 
is a homogeneous polynomial $p\in\K[x_{i_1},\dots,x_{i_r}]$ 
of degree $\deg (p)=\deg (m)$
satisfying the {\it generosity condition}
\begin{equation}
\K[x_{i_1},\dots,x_{i_r}]=\<R_{m}\>\oplus(p),
\label{eq:genericity}
\end{equation}
where $\<R_{m}\>$ is the linear span of the set $R_{m}$ of monomials 
in $\K[x_{i_1},\dots,x_{i_r}]$ which are not divisible by $m$, 
$(p)$ is the ideal in $\K[x_{i_1},\dots,x_{i_r}]$ generated by
$p$, and ``$\oplus$'' stands for a 
direct sum of subspaces. 
Notice that the generosity condition is satisfied for
a Zarisky open set of polynomials 
in $\K[x_{i_1},\dots,x_{i_r}]$ of degree $\deg(m)$.
For example, the polynomial $p=a\,x_1+b\,x_2$ is a $\{1,2\}$-deformation 
of the monomial $m=x_1$ if and only if $a\ne 0$.

The following lemma describes a class of $I$-deformations of monomials.

\begin {lemma}
Let $I=\{i_1,\dots,i_r\}$, let $m$ be a monomial in 
$\K[x_{i_1},\dots,x_{i_r}]$,
and let $\alpha_1,\dots,\alpha_r\in\K\setminus\{0\}$.
Then the polynomial
$$
p= (\alpha_1 x_{i_1}+\dots + \alpha_r x_{i_r})^{\deg m }
$$
is an $I$-deformation of the monomial $m$.
\label{le:generic}
\end{lemma}

\begin{proof}  
Let $m=x_{i_1}^{a_1}\cdots x_{i_r}^{a_r}$.
The generosity condition~(\ref{eq:genericity})
is equivalent to the condition that
the  operator
$$
\begin{array}{l}
A: \K[x_{i_1},\dots,x_{i_r}]\to \K[x_{i_1},\dots,x_{i_r}]\\[.05in]
A:f\mapsto (\partial/\partial x_{i_1})^{a_1}\cdots 
(\partial/\partial x_{i_r})^{a_r}(p\cdot f)
\end{array}
$$
has zero kernel. 
Let us change the coordinates to $y_1=x_{i_1}$,\dots, $y_{r-1}=x_{i_{r-1}}$,
$y_r=\alpha_1 x_{i_1}+\dots+\alpha_r x_{i_r}$. The operator $A$ can 
be written in these coordinates as
$$
A(f)=(\tilde\partial_1+\alpha_1\tilde\partial_r)^{a_1}\cdots
(\tilde \partial_{r-1}+\alpha_{r-1}\tilde \partial_r)^{a_{r-1}}
(\alpha_r\tilde \partial_r)^{a_r} (y_r^{a_1+\cdots + a_r}\cdot f)
$$
where $\tilde \partial_j=\partial/\partial y_j$. 
Then $A(f)=c\cdot f+g$, where $c$ is a nonzero constant and 
$\deg_{y_r}(g)<\deg_{y_r}(f)$.
Thus, in an appropriate basis, the operator $A$ is given by a 
triangular matrix with nonzero diagonal elements.
This implies that $\mathrm{Ker}\, A=0$.
\end{proof}

A {\it deformation} of a monotone monomial ideal $\I=\<m_I\mid I\in\Sigma\>$
is an ideal $\J=\<p_I\mid I\in\Sigma\>$ generated by polynomials $p_I$
such that $p_I$ is an $I$-deformation of $m_I$ for each
$I\in\Sigma$.  For example, according to Lemma~\ref{le:generic}, the ideal
$\J_G$ given in Section~\ref{sec:power} is a deformation of the monotone
monomial ideal $\I_G$.

\begin {theorem}
Let $\I$ be a monotone monomial ideal, 
and $R$ be the standard monomial basis of 
the algebra $\A=\K[x_1,\dots,x_n]/\I$, 
i.e., $R$ is the set of monomials that do not belong to $\I$.
Let $\J$ be a deformation of the ideal $\I$, and
$\B=\K[x_1,\dots,x_n]/\J$.

Then the monomials in $R$ linearly span the algebra $\B$.
\label{th:span}
\end{theorem}

Remark that the set of monomials $R$ may or may not be a basis for $\B$.

\begin{corollary}
\label{cor:I<J}  
Let $\I$ be a monotone monomial ideal, $\J$ 
be a deformation of the ideal $\I$,
$\A=\K[x_1,\dots,x_n]/\I$, and
$\B=\K[x_1,\dots,x_n]/\J$.
Then we have the 
following termwise inequalities for the Hilbert series:
$$
\Hilb\,\I\le \Hilb\,\J
\qquad\textrm{or equivalently,}\qquad
\Hilb\,\A\ge \Hilb\,\B.
$$ 
\end{corollary}

In some cases the Hilbert series are actually equal to each other. 
According to Theorem~\ref{th:ABK}, $\Hilb\, \A_G=\Hilb\,\B_G$, 
for any graph $G$.
However in general the Hilbert series may not be equal to each other.
It would be interesting to describe a general class of 
monotone monomial ideals and their deformations with equal
Hilbert series.

There is an obvious correspondence between the  generators $m_I$ of a 
monotone monomial ideal $\I$ and the generators $p_I$ of its 
deformation $\J$. 
Notice however that (except for very special cases)
the monomial generator $m_I$ does not belong to the boundary
of the Newton polytopes of it polynomial deformation $p_I$.
Thus the monomial $m_I$ is not the leading term of 
the polynomial $p_I$ for any term order.  
This shows that the above results cannot be tackled by the standard 
Gr\"obner bases technique.

\section{Syzygies of order monomial ideals}
\label{sec:min-free-res}

In this section we introduce a class of ideals that extends monotone
monomial ideals and construct free resolutions for these ideals.

Let $P$ be a finite partially ordered set, or {\it poset}.
Let $\M=\{m_u\mid u\in P\}$ be a collection of monomials in the polynomial
ring $\K[x_1,\dots,x_n]$ labelled by elements of the poset $P$.
Also let $M_u$ denote the set of all monomials divisible by $m_u$.
Let us say that $\M$ is an {\it order monomial family}, and the ideal
$\I=\<\M\>$ generated by the monomials $m_u$ is an 
{\it order monomial ideal}, if the following condition is satisfied:
\smallskip
\begin{itemize}
\item[(OM)] For any pair $u,v\in P$, there exists an upper bound $w\in P$
of $u$ and $v$
such that $M_u\cap M_v\subseteq M_w$,
i.e., $m_w$ divides $\lcm(m_u,m_v)$.
\end{itemize}
\smallskip
Here an upper bound means an element $w$ such that $w\geq u$ and $w\geq v$
in $P$. 
In particular, this condition implies that the poset $P$ has a unique
maximal element.

Every monotone monomial family is an order monomial family labelled 
by the set $P=\Sigma$ of subsets in $\{1,\dots,n\}$ partially 
ordered by inclusion.  Indeed, condition (MM3) is equivalent 
to condition~(OM).

Let $S=\K[x_1,\dots,x_n]$ be the polynomial ring.
For a non-negative $n$-vector $a=(a_1,\dots,a_n)\in\Z^n$, 
let $S(-a)=x^a \, S$
denote the free $\Z^n$-graded $S$-submodule in $S$ generated by the monomial
$x^a=x_1^{a_1}\cdots x_n^{a_n}$.  This submodule 
is isomorphic to $S$ with the $\Z^n$-grading shifted by the vector $a$.
If $a\geq b$ componentwise then $S(-a)$ is a submodule of $S(-b)$ 
and we will write $S(-a)\hookrightarrow S(-b)$ to denote the natural 
multidegree-preserving embedding of $S$-modules.

For an order monomial family $\M=\{m_u\mid u\in P\}$ and 
a subset $U\subseteq P$ of elements of the poset $P$,
let $m_U =\lcm(m_u\mid u\in U)$ be the least 
common multiple of the monomials $m_u$, $u\in U$.
We assume that $m_\emptyset = 1$.
Also let $a_U\in\Z^n$ be the exponent vector of the monomial
$m_U$.

Let us define the {\it homological order complex\/} $C_*(\M)$ for 
an order monomial ideal $\I=\<\M\>$ as the sequence of $\Z^n$-graded 
$S$-modules
\begin{equation}
\cdots \,{\overset{\partial_4}\longrightarrow}\,
C_3 \,{\overset{\partial_3}\longrightarrow}\,
C_2 \,{\overset{\partial_2}\longrightarrow}\, 
C_1 \,{\overset{\partial_1}\longrightarrow}\, 
C_0 = S \,{\longrightarrow}\, S/\I\to 0
\label{eq:CCC}
\end{equation}
whose $k$-th component is 
$$
C_k = \bigoplus_{u_1\lneq \cdots\lneq u_k} S(-a_{\{u_1,\dots,u_k\}}),
$$
where the direct sum is over strictly increasing $k$-chains 
$u_1\lneq \cdots\lneq u_k$ in the poset $P$.  The 
differential $\partial_k:C_k\to C_{k-1}$ is defined on the component 
$S(-a_{\{u_1,\dots,u_k\}})$ as the alternating sum 
$\partial_k = \sum_{i=1}^k (-1)^i\, E_i$ 
of the multidegree-preserving
embeddings $E_i: S(-a_{\{u_1,\dots,u_k\}})\hookrightarrow 
S(-a_{\{u_1,\dots, \widehat{u_i},\dots,u_k\}})$
of $S$-modules,
where $u_1,\dots, \widehat{u_i},\dots,u_k$ denotes the sequence
with skipped $i$-th element.

\begin{theorem}
Let $\M$ be an order monomial family.
The homological order complex $\C_*(\M)$ is a free resolution of 
the order monomial ideal $\I=\<\M\>$.

If $m_{\{u_1,\dots,u_k\}}\ne m_{\{u_1,\dots,\widehat{u_i},\dots,u_k\}}$,
for any increasing chain $u_1\lneq\cdots\lneq u_k$ in the poset $P$ and 
$i=1,\dots,k$, then the homological order complex $\C_*(\M)$ is a minimal free 
resolution of the order monomial ideal $\I=\<\M\>$.
\label{th:resolution1}
\end{theorem}

The above construction of the homological order complex $C_*(\M)$ 
is an instance of the general construction of cellular complexes for
monomial ideals due to Bayer and Sturmfels~\cite{BaSt}.
Their cellular complexes are associated with 
cell complexes,\footnote{``Cellular complexes'' should not be
confused with ``cell complexes.''
The former are homological complexes 
and the latter are geometrical complexes.}  
whose faces are marked by certain monomials.
In our case the cell complex is the {\it geometrical order complex\/} 
$\Delta=\Delta(P)$ of the poset $P$.  It is the simplical complex 
whose faces correspond to nonempty strictly increasing chains in $P$:  
$$
\Delta(P) = \{\{u_1, \cdots ,u_k\} \subseteq P \mid
u_1\lneq \cdots \lneq u_k,\ k\geq 1\}.
$$
For example, if $P$ consists of all nonempty subsets in 
$\{1,\dots,n\}$ ordered by inclusion then
$\Delta(P)$ is the barycentric subdivision of the $(n-1)$-dimensional simplex.
The face $F$ of $\Delta(P)$ given by an increasing chain 
$u_1\lneq \cdots\lneq u_k$ 
is marked by the monomial $m_F=m_{\{u_1,\dots,u_k\}}$.  
For a monomial $m$, let $\Delta_{\leq m}$ denote the subcomplex 
of $\Delta(P)$ formed by the faces $F$ whose mark $m_F$ divides $m$:
$$
\Delta_{\leq m} = \{F\in\Delta(P)\mid m_F \textrm{ divides } m\}.
$$
The faces of $\Delta(P)$ are partially ordered by containment of
closures.  More precisely, $F\geq F'$ if the increasing chain
for the face $F'$ is a subchain of the chain for $F$.

A result of Bayer and Sturmfels \cite[Proposition~1.2]{BaSt}
on cellular complexes implies the following statement. 

\begin{lemma} 
The complex $C_*(\M)$ is exact if and only if 
$\Delta_{\leq m}$ is acyclic over $\K$ for any monomial $m$.
If, in addition, $m_{F}\ne m_{F'}$, for any $F\gneq F'$,
then the complex $C_*(\M)$ is a minimal free resolution of the ideal 
$\I=\<\M\>$.
\label{lem:BaSt}
\end{lemma}

Actually, the subcomplex $\Delta_{\leq m}$ is not only acyclic but
also contractible.  This follows from the following result of Narushima.

\begin{lemma} {\rm\cite{Naru}} \
Let $M_u$, $u\in P$, be a finite collection of subsets 
in some set $\MF$, whose index set is a poset $P$. 
Assume that, for any $u,v\in P$, 
$M_u\cap M_v\subseteq M_w$
for some upper bound $w\in P$ of $u$ and $v$.
Then, for any $m\in\MF$, the subcomplex of the order complex $\Delta(P)$ 
of $P$ 
formed by the following collection of nonempty increasing chains in $P$
$$
\{\{u_1\lneq \cdots \lneq u_k\} \mid m\in M_{u_1}\cap \cdots \cap M_{u_k},
\ k\geq 1\}
$$
is contractible.
\label{lem:Narushima}
\end{lemma}

Theorem~\ref{th:resolution1} easily follow from 
Lemmas~\ref{lem:BaSt} and~\ref{lem:Narushima}.

\begin{proof}[Proof of Theorem~\ref{th:resolution1}]
Let $\M=\{m_u\mid u\in P\}$ be an order monomial family.
In view of Lemma~\ref{lem:BaSt} it is enough to show that 
the subcomplex $\Delta_{\leq m}$ is contractible for any 
monomial $m$.  According to~(OM), the conditions of 
Lemma~\ref{lem:Narushima} are satisfied,
where $\MF$ is the set of all monomials in $\K[x_1,\dots,x_n]$. 
For an increasing chain $u_1\lneq\cdots\lneq u_k$ in $P$,
the intersection $M_{u_1}\cap \cdots\cap M_{u_k}$ is the set of 
monomials divisible by $m_{\{u_1,\dots,u_k\}}$.  Thus the contractible complex
from Lemma~\ref{lem:Narushima} is exactly the complex $\Delta_{\leq m}$.
Lemma~\ref{lem:BaSt} implies that the homological order complex $C_*(\M)$ is
exact. 
\end{proof}

Let us now assume that $\M = \{m_I\mid I\in\Sigma\}$,
$m_I=\prod_{i\in I} x_i^{\nu_I(i)}$, is a monotone monomial family.
In this case, 
for a strictly increasing chain of subsets
$I_1\subsetneq \cdots \subsetneq I_k$, 
the least common multiple $m_{\{I_1,\dots,I_k\}}
=\lcm(m_{I_1},\dots, m_{I_k})$
is given by
\begin{equation}
m_{\{I_1,\dots,I_k\}}=\prod_{i_1\in I_1} x_{i_1}^{\nu_{I_1}(i_1)}\times
\prod_{i_2\in I_2\setminus I_1} x_{i_2}^{\nu_{I_2}(i_2)}\times
\cdots\times
\prod_{i_k\in I_k\setminus I_{k-1}} x_{i_k}^{\nu_{I_k}(i_k)}.
\label{eq:lcmm}
\end{equation}
In other words, 
the exponent vector $a_{\{I_1,\dots,I_k\}}=(a_1,\dots,a_n)\in\Z^n$
of the monomial $m_{\{I_1,\dots,I_k\}}$ is given by 
$$
a_i=\left\{\begin{array}{ll}
\nu_{I_r}(i) &\textrm{if } i\in I_r\setminus I_{r-1},\\[.1in]
0 &\textrm{if } i\not\in I_k,
\end{array}\right.
$$
where we assume that $I_0=\emptyset$.

Let us say that the monotone monomial family $\M$
and the corresponding monotone monomial ideal $\I=\<\M\>$
are {\it strictly monotone\/} if the following additional conditions hold:
\smallskip
\begin{itemize}
\item[(SM1)]  The ideal $\I$ is {\it minimally\/} generated by 
the set of monomials $\{m_I\mid I\in\Sigma\}$, i.e., there
are no two elements $I\ne J$ in $\Sigma$ such that $m_I$ divides $m_J$. 
\smallskip
\item[(SM2)] For any $I\subsetneq J\subsetneq K$ in $\Sigma$,
there exists $i\in J\setminus I$ such that $\nu_J(i)>\nu_K(i)$.
\end{itemize}
\smallskip

For example, a monotone monomial family such that 
the inequality in (MM2) is always strict and $\nu_I(i)>0$, for 
any $I\in\Sigma$ and $i\in I$, will be strictly monotone. 

Conditions~(SM1) and~(SM2) are equivalent to the statement that
$m_{\{I_1,\dots,I_k\}}\ne m_{\{I_1,\dots,\widehat{I_i},\dots,I_k\}}$
for any increasing chain $I_1\subsetneq\cdots \subsetneq I_k$ in $\Sigma$
and $i=1,\dots,k$.
Thus Theorem~\ref{th:resolution1} specializes to the following statement.

\begin{corollary}
Let $\M$ be a monotone monomial family.
Then the homological order complex $\C_*(\M)$ is a free resolution of 
the ideal $\I=\<\M\>$.
If $\M$ is a strictly monotone monomial family, then
$\C_*(\M)$ is a minimal free resolution of the strictly monotone ideal~$\I$.
\label{cor:resolution2}
\end{corollary}

Homological order complexes are related to Scarf complexes of 
generic monomial ideals.  
Let $\I=\<m_1,\dots,m_r\>$ be an arbitrary ideal in the polynomial
ring $S=\K[x_1,\dots,x_n]$ minimally generated 
by monomials $m_1,\dots,m_r$, and let $m_U=\lcm(m_u\mid u\in U)$
for $U\subseteq\{1,\dots,r\}$.
The {\it geometrical Scarf complex\/} of $\I$ was defined by
Bayer, Peeva, and Sturmfels~\cite{BPS}
as the following simplicial complex:
$$
\Delta_\I^{\textrm{scarf}} = 
\{U\subseteq\{1,\dots,r\} \mid m_U\ne m_V \textrm{ for all subsets } V\ne U,
\textrm{ and } |U|\geq 1 \}.
$$
The corresponding cellular complex is called the 
{\it homological Scarf complex}.

\begin{lemma}
Let $\I=\<m_u\mid u\in P\>$ be an order monomial ideal.
Then the geometrical Scarf complex $\Delta_\I^{\textrm{scarf}}$ 
is a subcomplex of the geometrical order complex $\Delta(P)$.
\label{lem:subcomplex}
\end{lemma}

\begin{proof}
Let $U$ be a subset of elements in $P$.  Suppose that $U$ 
contains two incomparable elements $u$ and $v$.  Let us pick an upper
bound $w$ of $u$ and $v$ provided by condition~(OM).  Let 
$V=U\cup\{w\}$ if $w\not\in U$, or
$V=U\setminus\{w\}$ if $w\in U$.
Then, according to~(OM), $m_U=m_V$.  Thus $V$ does
not belong to the geometrical Scarf complex.

This implies that, for any $U\in\Delta_\I^{\textrm{scarf}}$,
all elements of $U$ are comparable with each other, i.e.,
$U$ is an increasing chain in the poset $P$.
Thus $U\in\Delta(P)$.
\end{proof}

Let us say that a monomial $m$ {\it strictly divides\/} a monomial $m'$, if
$m$ divides $m'$ and $\deg_{x_i} (m'/m)\ne 0$ whenever $\deg_{x_i}(m')\ne 0$.
According to Miller, Sturmfels, and Yanagawa~\cite{MSY}, the ideal 
$\I=\<m_1,\dots,m_r\>$ is called {\it generic monomial ideal\/} 
if the following condition holds: 
\smallskip
\begin{itemize}
\item[(GM)] If two distinct minimal generators 
$m_u$ and $m_v$ have the same positive degree in some variable $x_i$, 
there is a third generator $m_w$ which strictly divides $\lcm(m_u,m_v)$.
\end{itemize}
\smallskip

The general, the Scarf complex may not be acyclic, but, for generic
monomial ideals, Miller, Sturmfels, and Yanagawa~\cite{MSY} 
proved the following result.

\begin{proposition} {\rm \cite[Corollary~1.8]{MSY}} \
If $\I$ is a generic monomial ideal, then the homological Scarf
complex is a minimal free resolution of $\I$.
\end{proposition}

We will see in Section~\ref{sec:three-examples}
that there are strictly monotone monomial ideals
that are not generic and there are generic monomial ideals that are
not strictly monotone.  The following claim
shows that these two classes of ideals have an interesting intersection.

\begin{proposition}  Let $\I=\<m_I\mid I\in\Sigma\>$ be a 
monotone monomial ideal such that 
the inequality in~{\rm(MM2)} is always strict, and $\nu_I(i)>0$, 
for any $I\in\Sigma$ and $i\in I$.
Then the monomial ideal $\I$ is both generic and strictly monotone.  
In this case, the geometrical/homological order complex for $\I$ 
coincides with geometrical/homological Scarf complex for $\I$.
\label{prop:order=scarf}
\end{proposition}

\begin{proof}
If monomials $m_I$ and $m_J$ have the same positive degree in some variable
then $I$ and $J$ are incomparable in $\Sigma$: 
$I\not\subseteq J$ and $J\not\subseteq I$.  By~(MM3) there exists
$K\supseteq I\cup J$ such that $m_K$ divides $\lcm(m_I,m_J)$.  
Then $K\ne I,J$. 
Since we assume that the inequality in~(MM2) is strict, 
$m_K$ strictly divides $\lcm(m_I,m_J)$.
It follows that the ideal $\I$ is generic.

According to Lemma~\ref{lem:subcomplex}, the geometrical Scarf complex
$\Delta_\I^{\textrm{scarf}}$ is a subcomplex of the geometrical order 
complex $\Delta(\Sigma)$.  Let us prove that 
$\Delta_\I^{\textrm{scarf}}=\Delta(\Sigma)$.
We need to show that, for any increasing chain
$I_1\subsetneq\cdots \subsetneq I_k$ in $\Sigma$, we have
$m_{\{I_1,\dots,I_k\}}\ne m_R$, where $R$ is any subset
of $\Sigma$ such that $R\ne\{I_1,\dots,I_k\}$.
This is clear if $R$ is a subchain in 
$I_1\subsetneq\cdots \subsetneq I_k$.
Otherwise, suppose that $m_{\{I_1,\dots,I_k\}}= m_R$ and
$R$ contains an element $J\not\in\{I_1,\dots,I_k\}$.
Then $m_J$ divides $m_{\{I_1,\dots,I_k\}}$.  According 
to conditions of the proposition, the monomial $m_J$ depends nontrivially 
on all $x_i$, $i\in J$.  Thus $J\subseteq I_k$.  Then 
$J\subseteq I_r$ and $J\not\subseteq I_{r-1}$ for some 
$r\in\{1,\dots,k\}$.  (We assume that $I_0=\emptyset$.)  
Pick an element $i\in J\setminus I_{r-1}$.  Then
$\nu_J(i) \leq \deg_{x_i} (m_{\{I_1,\dots,I_k\}}) = \nu_{I_r}(i)$
because $m_J$ divides $m_{\{I_1,\dots,I_k\}}$.
Also we have $\nu_J(i) > \nu_{I_r}(i)$ because $J\subsetneq I_r$.
Contradiction.
\end{proof}

The $k$-th {\it Betti number\/} $\beta_k(\I)$ of an ideal $\I$ is the rank of 
the $k$-th term in a minimal free resolution of $\I$.
The {\it graded Betti number\/} $\beta_{k,d}(\I)$ of a graded ideal
$\I$ is the number of direct summands in the $k$-th term of a 
minimal free resolution of $\I$ with generator of degree $d$.
Then $\beta_k(\I)=\sum_d \beta_{k,d}(\I)$.

Let $d(I_1,\dots,I_k)$ be the degree of the monomial $m_{\{I_1,\dots,I_k\}}$
given by
\begin{equation}
d(I_1,\dots,I_k) = \sum_{i_1\in I_1} \nu_{I_1}(i_1) +
\sum_{i_2\in I_2\setminus I_1} \nu_{I_2}(i_2)+\cdots+\
\sum_{i_k\in I_k\setminus I_{k-1}} \nu_{I_k}(i_k).
\label{eq:dIII}
\end{equation}

\begin{corollary}
Let $\I=\<m_I\mid I\in\Sigma\>$ be a strictly monotone monomial ideal.
The $k$-th Betti number $\beta_k(\I)$ of $\I$ is equal to the number of 
strictly increasing $k$-chains in the poset $\Sigma$.
Moreover, the graded Betti number $\beta_{k,d}(\I)$ of $\I$ is equal to the 
number of strictly increasing $k$-chains 
$I_1\subsetneq\cdots\subsetneq I_k$ in $\Sigma$
such that $d(I_1,\dots,I_k)=d$.

In particular, if $\Sigma$ is the set of all nonempty subsets in
$\{1,\dots,n\}$ then
\begin{equation}
\beta_k(\I) = k!\,S_{n+1,k+1},
\label{eq:Betti-Stirling}
\end{equation}
where $S_{n+1,k+1}$ is the Stirling number of the second kind,
i.e., the number of partitions of the set $\{0,\dots,n\}$
into $k+1$ nonempty blocks.
\label{cor:graded-Betti-Stirling}
\end{corollary}

The last claim is obtained by counting strictly increasing $k$-chains 
of nonempty subsets $I_1\subsetneq \cdots \subsetneq I_k$ 
in $\{1,\dots,n\}$.  Indeed, 
these chains are in one-to-one correspondence with 
partitions of the set $\{0,\dots,n\}$ into a linearly  ordered family of $k+1$ 
nonempty blocks $(I_1, I_2\setminus I_1,\cdots, I_{k}\setminus I_{k-1}, 
\{0,\dots,n\}\setminus I_k)$ such that the last block contains $0$.
There are $k!$ ways to pick such a linear ordering of blocks.

Let us say that a (directed) graph is {\it saturated\/} if all off-diagonal 
entries of the adjacency matrix $A=(a_{ij})$ are nonzero: 
$a_{ij}\geq 1$ for $i\ne j$.  
If $G$ is a saturated digraph then the monotone monomial ideal 
$\I_G$ constructed in Section~\ref{sec:G-parking} 
satisfies the conditions of Proposition~\ref{prop:order=scarf}.

\begin{corollary}
The monotone monomial ideal $\I_G$, for a saturated digraph $G$, is both
strictly monotone and generic.  In this case $\Sigma$ is the poset 
of all nonempty subsets in $\{1,\dots,n\}$.  The homological order complex 
$C_*(\M)$, which coincides with the homological Scarf complex, gives a minimal 
free resolution of the ideal $\I_G$. 
Its Betti numbers are given by formula~{\rm(\ref{eq:Betti-Stirling})}.
\label{cor:IG-saturated}
\end{corollary}

It would be interesting to find a minimal free resolution of the ideal 
$\I_G$ for any non-saturated digraph $G$.  More generally, it would be
interesting to find a minimal free resolution for any monotone monomial
ideal.

Computer experiments suggest the following conjecture on Betti numbers 
of deformations of monotone monomial ideals.

\begin{conjecture}  Let $\J$ be a deformation of a monotone
monomial ideal $\I$ such that 
$\dim_{\K} \K[x_1,\dots,x_n]/\I = \dim_{\K} \K[x_1,\dots,x_n]/\J$.
Then all graded Betti numbers of the ideals $\I$ and $\J$ coincide:
$\beta_{k,d}(\I) = \beta_{k,d}(\J)$.
In particular, for a graph $G$,
the ideals $\I_G$ and $\J_G$ have the same graded Betti numbers.
\label{conj:Betti=Betti}
\end{conjecture}

\section{Examples}
\label{sec:three-examples}

Let us illustrate Corollaries~\ref{cor:resolution2}, 
\ref{cor:graded-Betti-Stirling}, and~\ref{cor:IG-saturated}
and Proposition~\ref{prop:order=scarf}
by several examples.
In all examples $n=3$,  $S=\K[x_1,x_2,x_3]$, and 
$S(-d)$ denotes the $\Z$-graded $S$-module isomorphic to $S$ 
with grading shifted by integer $d$, so that the generator
has degree~$d$.

\begin{example} {\rm 
Let $G=K_4$ be the complete graph on 4 vertices.  This graph is saturated.  
Thus the monomial ideal $\I_G$ is both strictly monotone and generic.
The poset $\Sigma$ consists of all nonempty subsets in $\{1,2,3\}$.  
The Hasse diagram of $\Sigma$ marked by the monomials
$m_I$ is given by 
$$
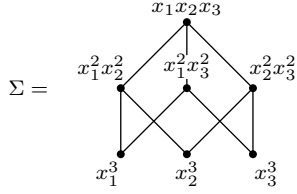
\pspicture(0,-25)(100,120)
\pscircle*(0,0){3}
\pscircle*(50,0){3}
\pscircle*(100,0){3}
\pscircle*(0,50){3}
\pscircle*(50,50){3}
\pscircle*(100,50){3}
\pscircle*(50,100){3}
\psline{-}(0,0)(0,50)
\psline{-}(0,0)(50,50)
\psline{-}(50,0)(0,50)
\psline{-}(50,0)(100,50)
\psline{-}(100,0)(100,50)
\psline{-}(100,0)(50,50)
\psline{-}(0,50)(50,100)
\psline{-}(50,50)(50,57)
\psline{-}(50,75)(50,100)
\psline{-}(100,50)(50,100)

\rput(-70,50){{\footnotesize $\Sigma=$}}
\rput(-10,-15){{\footnotesize $x_1^3$}}
\rput(50,-15){{\footnotesize $x_2^3$}}
\rput(110,-15){{\footnotesize $x_3^3$}}
\rput(-15,65){{\footnotesize $x_1^2 x_2^2$}}
\rput(115,65){{\footnotesize $x_2^2 x_3^2$}}
\rput(50,67){{\footnotesize $x_1^2 x_3^2$}}
\rput(50,110){{\footnotesize $x_1 x_2 x_3$}}
\endpspicture
$$
The poset $\Sigma$ has seven 1-chains, twelve 2-chains, and six 3-chains.
In this case the geometrical order complex $\Delta=\Delta(\Sigma)$ 
is the barycentric subdivision of a triangle.  
The following figure shows the complex $\Delta$
with faces marked by vectors $a_{\{I_1,\dots,I_k\}}$:
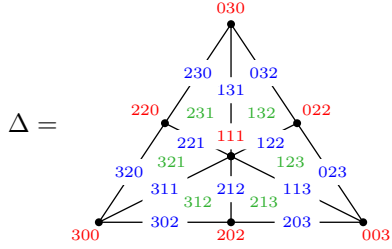
$$
\pspicture(0,-20)(200,170)
\pscircle*(0,0){3}
\pscircle*(100,0){3}
\pscircle*(200,0){3}
\pscircle*(100,150){3}
\pscircle*(100,50){3}
\pscircle*(100,50){3}
\pscircle*(50,75){3}
\pscircle*(150,75){3}
\psline{-}(0,0)(200,0)(100,150)(0,0)
\psline{-}(0,0)(150,75)
\psline{-}(200,0)(50,75)
\psline{-}(100,0)(100,150)
\rput(-10,-10){\red{\tiny$300$}}
\rput(210,-10){\red{\tiny$003$}}
\rput(100,162){\red{\tiny$030$}}
\rput(100,-10){\red{\tiny$202$}}
\rput(35,85){\red{\tiny$220$}}
\rput(165,85){\red{\tiny$022$}}
\rput(100,65){\red{\psframebox*{\tiny$111$}}}
\rput(100,25){\blue{\psframebox*{\tiny$212$}}}
\rput(50,25){\blue{\psframebox*{\tiny$311$}}}
\rput(150,25){\blue{\psframebox*{\tiny$113$}}}
\rput(50,0){\blue{\psframebox*{\tiny$302$}}}
\rput(150,0){\blue{\psframebox*{\tiny$203$}}}
\rput(100,100){\blue{\psframebox*{\tiny$131$}}}
\rput(23,38){\blue{\psframebox*{\tiny$320$}}}
\rput(177,38){\blue{\psframebox*{\tiny$023$}}}
\rput(75,113){\blue{\psframebox*{\tiny$230$}}}
\rput(125,113){\blue{\psframebox*{\tiny$032$}}}
\rput(70,60){\blue{\psframebox*{\tiny$221$}}}
\rput(130,60){\blue{\psframebox*{\tiny$122$}}}
\rput(55,46){\mygreen{\tiny$321$}}
\rput(145,46){\mygreen{\tiny$123$}}
\rput(75,15){\mygreen{\tiny$312$}}
\rput(125,15){\mygreen{\tiny$213$}}
\rput(77,83){\mygreen{\tiny$231$}}
\rput(123,83){\mygreen{\tiny$132$}}
\rput(-50,75){$\Delta=$}
\endpspicture
$$
The Betti numbers $(\beta_0,\beta_1,\beta_2,\beta_3) = (1,7,12,6)$
of the ideal $\I=\I_{K_4}$,
which are also the $f$-numbers of the order complex $\Delta$,
can be expressed in terms of the Stirling numbers 
by formula~(\ref{eq:Betti-Stirling}).
The graded Betti numbers of this ideal are indicated
on the following minimal free resolution:
$$
0
\longrightarrow S(-6)^{6}
\longrightarrow S(-5)^{12}
\longrightarrow S(-3)^{4}\oplus S(-4)^{3}
\longrightarrow S 
\longrightarrow S/\I 
\longrightarrow 0.
$$
This resolution is the homological order complex and also 
the homological Scarf complex of~$\I$.

Similarly, 
a minimal free resolution of the ideal $\I_n=\I_{K_{n+1}}$ 
associated with the complete graph $K_{n+1}$ is given by 
the cellular complex corresponding to the simplicial complex
$\Delta(\Sigma)=\Delta^{\textrm{scarf}}_{\I_n}$, which is  
the barycentric subdivision of the $(n-1)$-dimensional simplex,
cf.~\cite[Example~1.2]{MSY}.
}
\label{ex:res1}
\end{example}

\begin{example} {\rm
Let $G$ be the graph given by
$$
G=
\lower.2in
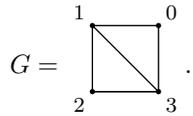
\hbox{\pspicture(-20,-15)(70,67)
\psline{-}(0,0)(50,0)(50,50)(0,50)(0,0)
\psline{-}(0,0)(50,0)(0,50)
\pscircle*(0,0){2}
\pscircle*(50,0){2}
\pscircle*(0,50){2}
\pscircle*(50,50){2}
\rput(60,60){{\footnotesize0}}
\rput(-10,60){{\footnotesize1}}
\rput(-10,-10){{\footnotesize2}}
\rput(60,-10){{\footnotesize3}}
\endpspicture}.
$$
This graph is not saturated and the monotone monomial family that 
generates the ideal $\I=\I_G
=\<x_1^3,x_2^2,x_3^3,x_1^2x_2,x_1^2x_3^2,x_2x_3^2,x_1x_2^0x_3\>$
will not be strictly monotone if we 
assume that the labelling set $\Sigma$ consists 
of all nonempty subsets in $\{1,2,3\}$.
As we mentioned before, the generator $m_{\{1,3\}}=x_1^2 x_3^2$
is redundant.  The same monomial ideal $\I$ 
is minimally generated by the {\it strictly\/} 
monotone monomial family $\{m_I\mid I\in\Sigma\}$ with $\Sigma=\{\{1\},
\{2\},\{3\},\{1,2\},\{2,3\},\{1,2,3\}\}$.
The Hasse diagram of this poset $\Sigma$ marked 
by the monomials $m_I$ is given by
$$
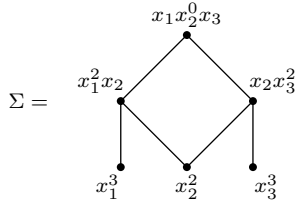
\pspicture(0,-25)(100,125)
\pscircle*(0,0){3}
\pscircle*(50,0){3}
\pscircle*(100,0){3}
\pscircle*(0,50){3}
\pscircle*(100,50){3}
\pscircle*(50,100){3}
\psline{-}(0,0)(0,50)
\psline{-}(50,0)(0,50)
\psline{-}(50,0)(100,50)
\psline{-}(100,0)(100,50)
\psline{-}(0,50)(50,100)
\psline{-}(100,50)(50,100)

\rput(-70,50){{\footnotesize $\Sigma=$}}
\rput(-10,-15){{\footnotesize $x_1^3$}}
\rput(50,-15){{\footnotesize $x_2^2$}}
\rput(110,-15){{\footnotesize $x_3^3$}}
\rput(-15,65){{\footnotesize $x_1^2 x_2$}}
\rput(115,65){{\footnotesize $x_2 x_3^2$}}
\rput(50,115){{\footnotesize $x_1 x_2^0 x_3$}}
\endpspicture
$$
This poset has six 1-chains, nine 2-chains, and four 3-chains.
Its geometrical order complex $\Delta=\Delta(\Sigma)$ with faces marked 
by vectors $a_{I_1,\dots,I_k}$
is shown on the following figure:
$$
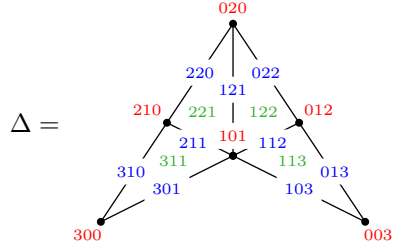
\pspicture(0,-20)(200,170)
\pscircle*(0,0){3}
\pscircle*(200,0){3}
\pscircle*(100,150){3}
\pscircle*(100,50){3}
\pscircle*(100,50){3}
\pscircle*(50,75){3}
\pscircle*(150,75){3}
\psline{-}(200,0)(100,150)(0,0)
\psline{-}(0,0)(150,75)
\psline{-}(200,0)(50,75)
\psline{-}(100,50)(100,150)
\rput(-10,-10){\red{\tiny$300$}}
\rput(210,-10){\red{\tiny$003$}}
\rput(100,162){\red{\tiny$020$}}
\rput(35,85){\red{\tiny$210$}}
\rput(165,85){\red{\tiny$012$}}
\rput(100,65){\red{\psframebox*{\tiny$101$}}}
\rput(50,25){\blue{\psframebox*{\tiny$301$}}}
\rput(150,25){\blue{\psframebox*{\tiny$103$}}}
\rput(100,100){\blue{\psframebox*{\tiny$121$}}}
\rput(23,38){\blue{\psframebox*{\tiny$310$}}}
\rput(177,38){\blue{\psframebox*{\tiny$013$}}}
\rput(75,113){\blue{\psframebox*{\tiny$220$}}}
\rput(125,113){\blue{\psframebox*{\tiny$022$}}}
\rput(70,60){\blue{\psframebox*{\tiny$211$}}}
\rput(130,60){\blue{\psframebox*{\tiny$112$}}}
\rput(55,46){\mygreen{\tiny$311$}}
\rput(145,46){\mygreen{\tiny$113$}}
\rput(77,83){\mygreen{\tiny$221$}}
\rput(123,83){\mygreen{\tiny$122$}}
\rput(-50,75){$\Delta=$}
\endpspicture
$$
The homological order complex give a minimal free resolution
of the ideal $\I$:
$$
0
\longrightarrow S(-5)^{4}
\longrightarrow S(-4)^{9}
\longrightarrow S(-2)^{2}\oplus S(-3)^{4}
\longrightarrow S 
\longrightarrow S/\I 
\longrightarrow 0.
$$
In this case the ideal $\I$ is also generic and the above resolution
is the homological Scarf complex.
}
\label{ex:res2}
\end{example}

\begin{example} {\rm
Let $\I=\<x_1^2x_2^2,x_2^2x_3,x_1x_2x_3\>$ be the monotone monomial
ideal, whose poset $\Sigma$ marked by the monomials is given by
$$
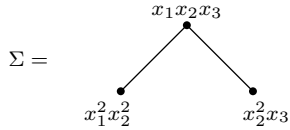
\pspicture(0,-25)(100,70)
\pscircle*(0,0){3}
\pscircle*(100,0){3}
\pscircle*(50,50){3}
\psline{-}(0,0)(50,50)(100,0)
\rput(-70,25){{\footnotesize $\Sigma=$}}
\rput(-10,-17){{\footnotesize $x_1^2x_2^2$}}
\rput(110,-17){{\footnotesize $x_2^2x_3$}}
\rput(50,60){{\footnotesize $x_1x_2x_3$}}
\endpspicture
$$
The ideal $\I$ is strictly monotone but {\it is not\/} generic.
The geometrical order complex with faces marked by vectors
$a_{\{I_1,\dots,I_k\}}$ is given by
$$
\pspicture(0,-10)(200,20)
\pscircle*(0,0){3}
\pscircle*(100,0){3}
\pscircle*(200,0){3}
\psline{-}(0,0)(200,0)
\rput(0,10){\red{\tiny$220$}}
\rput(100,10){\red{\tiny$111$}}
\rput(200,10){\red{\tiny$021$}}
\rput(50,10){\blue{\tiny$221$}}
\rput(150,10){\blue{\tiny$121$}}
\rput(-50,0){$\Delta=$}
\endpspicture
$$
It produces the following minimal free resolution of the ideal $\I$:
$$
0
\longrightarrow S(-4)\oplus S(-5)
\longrightarrow S(-3)^{2}\oplus S(-4)
\longrightarrow S 
\longrightarrow S/\I 
\longrightarrow 0.
$$
On the other hand, the geometrical Scarf complex in this case is disconnected
and does not give a resolution for $\I$.
}
\label{ex:res3}
\end{example}

\begin{example} {\rm
Let $G$ be the graph given by
$$
G=
\lower.2in
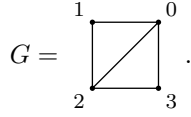
\hbox{\pspicture(-20,-10)(70,67)
\psline{-}(0,0)(50,0)(50,50)(0,50)(0,0)
\psline{-}(0,0)(0,0)(50,50)
\pscircle*(0,0){2}
\pscircle*(50,0){2}
\pscircle*(0,50){2}
\pscircle*(50,50){2}
\rput(60,60){{\footnotesize0}}
\rput(-10,60){{\footnotesize1}}
\rput(-10,-10){{\footnotesize2}}
\rput(60,-10){{\footnotesize3}}
\endpspicture}.
$$
Then $\I=\I_G 
= \<x_1^2,x_2^3,x_3^2,x_1 x_2^2,x_1^2 x_3,x_1^2 x_3^2,x_1 x_2 x_3\>$.
The Hasse diagram of $\Sigma$ 
marked by the monomials $m_I$ is given by
$$
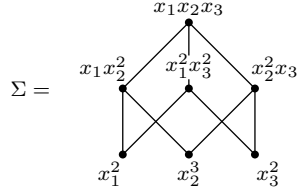
\pspicture(0,-25)(100,120)
\pscircle*(0,0){3}
\pscircle*(50,0){3}
\pscircle*(100,0){3}
\pscircle*(0,50){3}
\pscircle*(50,50){3}
\pscircle*(100,50){3}
\pscircle*(50,100){3}
\psline{-}(0,0)(0,50)
\psline{-}(0,0)(50,50)
\psline{-}(50,0)(0,50)
\psline{-}(50,0)(100,50)
\psline{-}(100,0)(100,50)
\psline{-}(100,0)(50,50)
\psline{-}(0,50)(50,100)
\psline{-}(50,50)(50,57)
\psline{-}(50,75)(50,100)
\psline{-}(100,50)(50,100)

\rput(-70,50){{\footnotesize $\Sigma=$}}
\rput(-10,-15){{\footnotesize $x_1^2$}}
\rput(50,-15){{\footnotesize $x_2^3$}}
\rput(110,-15){{\footnotesize $x_3^2$}}
\rput(-15,65){{\footnotesize $x_1 x_2^2$}}
\rput(115,65){{\footnotesize $x_2^2 x_3$}}
\rput(50,67){{\footnotesize $x_1^2 x_3^2$}}
\rput(50,110){{\footnotesize $x_1 x_2 x_3$}}
\endpspicture
$$
The corresponding homological order complex $C_*(\M)$ gives 
a free resolution, which is not minimal.  In this case 
the monomial generator $m_{\{1,3\}} = x_1^2 x_3^2$ is redundant
because it is divisible by $m_{\{1\}} = x_1^2$.
However, the family $\M\setminus \{m_{\{1,3\}}\}$ {\it is not\/}
a monotone monomial family because condition (MM3) 
no longer holds.  
On the other hand, $\I$ is a generic ideal and its Scarf complex
gives a minimal free resolution:
$$
0
\longrightarrow S(-5)^{4}
\longrightarrow S(-4)^{9}
\longrightarrow S(-2)^{2}\oplus S(-3)^{4}
\longrightarrow S 
\longrightarrow S/\I 
\longrightarrow 0.
$$
}
\label{ex:res4}
\end{example}

In Examples~\ref{ex:res1}, \ref{ex:res2}, and~\ref{ex:res4} above, the 
graded Betti numbers of the deformed ideal $\J_G$ coincide 
with graded Betti numbers of $\I_G$.

\section{Hilbert series and dimensions of monotone monomial ideals}
\label{sec:hilbert-monomial}

In this section we give formulas for the Hilbert series 
and dimensions of monotone monomial ideals.  Then we prove
Theorem~\ref{th:G-parking=trees}.

Let $\{m_I\mid I\in\Sigma\}$, $m_I=\prod_{i\in I} x_i^{\nu_I(i)}$, 
be a monotone monomial family, let $\I=\<m_I\>$ be the corresponding ideal in 
the polynomial ring $\K[x_1,\dots,x_n]$, and let $\A=\K[x_1,\dots,x_n]/\I$.

\begin{proposition}  The Hilbert series of the algebra $\A$ equals 
$$
\Hilb\,\A = 
\frac{\ds 1+\sum_{k\geq 1}(-1)^k\sum_{I_1\subsetneq \cdots \subsetneq I_k}
q^{d(I_1,\dots,I_k)}}
{(1-q)^n}\,,
$$
where the sum is over all strictly increasing chains 
in $\Sigma$ and the number $d(I_1,\dots,I_k)$ is given 
by~{\rm(\ref{eq:dIII})}.
\label{prop:monotone-GF}
\end{proposition}

\begin{proof}[First Proof of Proposition~\ref{prop:monotone-GF}]
According to Corollary~\ref{cor:resolution2}, the 
homological order complex~(\ref{eq:CCC}) is a free resolution of $\I$.
Thus 
$$
\Hilb\,\A = \sum_{k\geq 0} (-1)^k \,\Hilb\,C_k
=\sum_{k\geq 0} (-1)^k \sum_{I_1\subsetneq \cdots \subsetneq I_k} 
\Hilb\,S(-a_{\{I_1,\dots,I_k\}}).
$$
Since $\Hilb\,S(-a_{\{I_1,\dots,I_k\}}) = q^{d(I_1,\dots,I_k)}/(1-q)^n$,
the proposition follows.
\end{proof}

Let us give another more expanded proof of Proposition~\ref{prop:monotone-GF}.
We will need the {\it improved inclusion-exclusion formula\/} 
due to Narushima~\cite{Naru}.
For a subset $M$ in some set $\MF$, 
let $\chi(M)$ denote the characteristic function of $M$:
$$
\chi(M)\,:\,a\longmapsto\left\{\begin{array}{ll}
1 &\textrm{if } a\in M,\\[.1in]
0 &\textrm{if } a\in \MF\setminus M.
\end{array}\right.
$$

\begin{lemma} {\rm \cite{Naru}} \
Let $M_u$, $u\in P$, be a finite collection of subsets 
in some set $\MF$, whose index set is a poset $P$, 
such that, for any $u,v\in P$, $M_u\cap M_v \subseteq M_w$ for some upper
bound $w$ of $u$ and $v$.
Then we have
$$
\chi(\MF\setminus\bigcup_{u\in P} M_u) =
\chi(\MF) + \sum_{k\geq 1}(-1)^k \sum_{u_1<\cdots<u_k} 
\chi(M_{u_1}\cap \cdots \cap M_{u_k}),
$$
where the second sum is over all strictly increasing chains
$u_1<\cdots<u_k$ in the poset $P$.
\label{lem:GFMX}
\end{lemma}

\begin{proof}
According to the usual {\it inclusion-exclusion principle,} 
see~\cite[Section~2.1]{EC1}, we have
$$
\chi(\MF\setminus\bigcup_{u\in P} M_u) = 
\chi(\MF) - \sum_u \chi(M_u) + \sum_{u,v} \chi(M_u\cap M_v) - \cdots.
$$
The general summand in this expression is
$s_{U} = (-1)^k\, \chi(M_{u_1}\cap \cdots \cap M_{u_k})$,
where $U=\{u_1,\dots,u_k\}$ is an unordered $k$-element subset in $P$.
We argue that if we take the summation only over {\it
increasing chains\/} $u_1<\cdots < u_k$ in $P$
we get exactly the same answer.
Indeed, let us show that the contribution of all other subsets 
$U$ is zero.  We will use the {\it involution principle,}
see~\cite[Section~2.6]{EC1}.
Let us construct an involution $\iota$ on the set of all subsets 
$U\subseteq P$ of all possible sizes $k\geq 0$ such that 
the elements of $U$ cannot be arranged in an increasing chain.
Let us fix a linear order on elements of $P$.
Find the lexicographically minimal pair of of incomparable elements 
$u$ and $v$ in $U$, i.e., $u\not\leq v$ and $v\not\leq u$.  
Let $w\in P$ be the minimal (with respect to the linear order) 
upper bound of $u$ and $v$ such that $M_u\cap M_v \subseteq M_w$.  
Define the map $\iota$ as follows:
$$
\iota\,:\,U\longmapsto
\left\{\begin{array}{ll}
 U\cup\{w\},&\textrm{if } w\not\in U,\\[.1in]
 U\setminus \{w\},&\textrm{if } w\in U.
\end{array}\right.
$$
Then $\iota$ is an involution such that $|\iota(U)|=|U|\pm 1$.
Conditions of the lemma imply that
$s_{U}=-s_{\iota(U)}$.  Thus all summands $s_{U}$ 
corresponding to non-chains cancel each other.
\end{proof}

\begin{proof}[Second Proof of Proposition~\ref{prop:monotone-GF}]
Let $\MF$ be the set of monomials in $\K[x_1,\dots,x_n]$,
and, for $I\in \Sigma$, let $M_I\subset \MF$ denote the set of monomials in 
$\K[x_1,\dots,x_n]$ divisible by $m_I$.
For a subset of monomials $M\subset\MF$, let
$$
[M] = \sum_{m\in M} q^{\deg(m)} = \sum_{m\in\MF}  q^{\deg(m)} \chi(M)(m).
$$
Then $\Hilb\,\A=[\MF\setminus\bigcup M_I]$ and
$[\MF]=(1-q)^{-n}$.
All conditions of Lemma~\ref{lem:GFMX} are satisfied,
where $P=\Sigma$.   
For an increasing chain $I_1\subsetneq \cdots \subsetneq I_k$, 
the least common multiple~(\ref{eq:lcmm}) of the monomials
$m_{I_1},\dots, m_{I_k}$ has degree $d(I_1,\dots,I_k)$.
Thus
$[M_{I_1}\cap \cdots \cap M_{I_k}] = q^{d(I_1,\dots,I_k)} (1-q)^{-n}$.
Proposition~\ref{prop:monotone-GF} follows from Lemma~\ref{lem:GFMX}.
\end{proof}

\begin{lemma}
The algebra $\A$ is finite-dimensional as a linear space over $\K$
if and only if $\Sigma$ contains all one-element subsets in 
$\{1,\dots,n\}$.
\end{lemma}

\begin{proof}
If there is $i\in\{1,\dots,n\}$ such that $\{i\}\not\in\Sigma$
then the powers $x_i^c$ form an infinite linearly 
independent subset in $\A$.  Thus $\A$ is infinite-dimensional.
Otherwise, if $\Sigma$ contains all one-element subsets, then
the algebra $\A$ is finite-dimensional.  Indeed, a monomial 
$x_1^{a_1}\dots x_n^{a_n}$ vanishes in $\A$
unless $a_1< \nu_{\{1\}}(1),\dots,a_n<\nu_{\{n\}}(n)$.
\end{proof}

\begin{proposition}
Assume that $\Sigma$ contains all one-element subsets
in $\{1,\dots,n\}$, and  let $\nu(i) = \nu_{\{i\}}(i)$.
The dimension of the algebra $\A$ is given by the following
polynomial in the variables $\{\nu_I(i)\mid I\in\Sigma, i\in I\}$:
\begin{equation}
\begin{array}{l}
\ds\dim \A = \sum_{I_1\subsetneq \cdots \subsetneq I_k} (-1)^k
\prod_{i_1\in I_1}(\nu(i_1)-\nu_{I_1}(i_1))\times
\prod_{i_2\in I_2\setminus I_1}(\nu(i_2)-\nu_{I_2}(i_2))\times
\\[.3in]
\ds\qquad\qquad\qquad\qquad
\times \cdots \times
\prod_{i_k\in I_k\setminus I_{k-1}}(\nu(i_k)-\nu_{I_k}(i_k))\times
\prod_{i_{k+1}\not\in I_k}\nu(i_{k+1}),
\end{array}
\label{eq:monotone-chains}
\end{equation}
where the sum is over all strictly increasing chains 
$I_1\subsetneq \cdots \subsetneq I_k$ of nonempty subsets in 
$\{1,\dots,n\}$ of all sizes $k\geq 0$, 
including the empty chain of size $k=0$.
\label{prop:monotone-chains}
\end{proposition}

\begin{proof}
Let $\tilde\MF$ be the set of monomials $x_1^{a_1}\cdots x_n^{a_n}$ 
such that $a_1< \nu(1),\dots,a_n<\nu(n)$.
A monomial $x^a$ vanishes in the algebra $\A$ unless $x^a\in\tilde\MF$.
Let $\tilde M_I = M_I\cap \tilde\MF$.
Lemma~\ref{lem:GFMX} for $\tilde\MF$ implies that
$$
\dim\A = 
|\tilde\MF|+\sum_{k\geq 1}(-1)^k 
\sum_{I_1\subsetneq \cdots \subsetneq I_k} 
|\tilde M_{I_1}\cap\cdots \cap \tilde M_{I_k}|.
$$
The intersection $\tilde M_{I_1}\cap\cdots \cap \tilde M_{I_k}$ is the set 
of all monomials in $\tilde\MF$ divisible by the monomial 
$m_{I_1,\dots,I_k}$ given by~(\ref{eq:lcmm}).
Thus $(-1)^k\,|\tilde M_{I_1}\cap\cdots \cap \tilde M_{I_k}|$
is equal to the summand in~(\ref{eq:monotone-chains}).
\end{proof}

Remark that if
$I_1,\dots,I_k$ is not a chain of subsets then
$|\tilde M_{I_1}\cap\cdots \cap \tilde M_{I_k}|$
may not be a polynomial in the $\nu_I(i)$.  It
may include expressions like $\min(\nu_I(i),\nu_J(i))$.
Thus the inclusion-exclusion principle does not immediately
produce a polynomial expression for $\dim\A$.
Miraculously, all non-polynomial terms cancel each other.

\medskip
We can now prove Theorem~\ref{th:G-parking=trees} that 
claims that dimension of the algebra $\A_G$ equals
the number of oriented spanning trees of $G$.

\begin{proof}[Proof of Theorem~\ref{th:G-parking=trees}]
cf.~\cite[Appendix~E]{Gab1} \
Let $G$ be a digraph on the vertices $0,\dots,n$, and let 
$A=(a_{ij})$ be its adjacency matrix.
Specializing Proposition~\ref{prop:monotone-chains}, we 
obtain the following polynomial
formula for the dimension of the algebra $\A_G$:
\begin{equation}
\begin{array}{l}
\ds\dim \A_G = \sum_{I_1\subsetneq \cdots \subsetneq I_k} (-1)^k
\prod_{i_1\in I_1}\(\sum_{j_1\in I_1} a_{i_1 j_1}\)\times
\prod_{i_2\in I_2\setminus I_1}\(\sum_{j_2 \in I_2} a_{i_2 j_2}\)\times
\\[.3in]
\ds
\times \cdots \times
\prod_{i_k\in I_k\setminus I_{k-1}}\(\sum_{j_k\in I_k} a_{i_k j_k}\)
\times
\prod_{i_{k+1}\in \{1,\dots,n\}\setminus I_k}\(\sum_{0\leq j_{k+1}\leq n}
a_{i_{k+1} j_{k+1}}\),
\end{array}
\label{eq:graph-chains}
\end{equation}
where the sum is over all strictly increasing chains 
$I_1\subsetneq \cdots \subsetneq I_k$ of nonempty subsets 
in $\{1,\dots,n\}$ of all sizes $k\geq 0$. 
In this formula, we assume that $a_{ii}=0$.

Let us show that the expression~(\ref{eq:graph-chains}) for $\dim \A_G$ 
is equal to the number of oriented spanning trees of $G$.
We will use the involution principle again.

Let us first give a combinatorial interpretation of the 
right-hand side of~(\ref{eq:graph-chains}).
The summand that corresponds to an increasing chain
$I_*=I_1\subsetneq \cdots \subsetneq I_k$ is equal to $(-1)^k$ times
the number of subgraphs $H$ of $G$ such that 
\begin{enumerate}
\item[(a)] $H$ contains exactly $n$ directed 
edges $(i,f(i))$ for $i=1,\dots,n$.
\item[(b)] If $i\in I_r\setminus I_{r-1}$ then $f(i)\in I_r$.
(We assume that $I_0=\emptyset$.)
\end{enumerate}

For such a subgraph $H$, let $J_H\subseteq\{1,\dots,n\}$ be the set of 
vertices $i$ such that $f^p(i)=0$ for some power $p$, i.e., $J_H$ 
is the set of vertices $i$ such that there is a directed path in $H$ 
from $i$ to the root $0$.
Notice that if $i\in\bigcup I_r$ then $f^p(i)\in\bigcup I_r$, 
thus $f^p(i)\ne 0$, for any $p$.
Thus $I_1,\dots,I_k\subseteq \overline{J_H}=\{1,\dots,n\}\setminus J_H$.
Also notice that $H$ is an oriented spanning tree of $G$ if and only if 
$J_H=\{1,\dots,n\}$.  

Let us now construct an involution $\kappa$ on the set of pairs
$(I_*, H)$ such that $H$ is not an oriented spanning tree.  
In other words, the involution $\kappa$ acts
on the set of pairs $(I_*, H)$ with nonempty $\overline{J_H}$.
It is given by
$$
\kappa\,:\,(I_*,H)\longmapsto\left\{
\begin{array}{ll}
(I_1\subsetneq \cdots \subsetneq I_{k-1}, H)
&\textrm{if } I_k= \overline{J_H}\,,\\[.1in]
(I_1\subsetneq \cdots \subsetneq I_k\subsetneq \overline{J_H}, H)
&\textrm{if } I_k\ne \overline{J_H}\,.
\end{array}\right.
$$
The contribution of the pair $(I_*, H)$ to the 
right-hand side of~(\ref{eq:graph-chains}) 
is opposite to the contribution of $\kappa(I_*, H)$. 
Thus the contributions of all subgraphs $H$ which are not oriented
spanning trees cancel each other.  This implies that $\dim \A_G$ 
is the number of oriented spanning trees.
\end{proof}

It would be interesting to find a combinatorial proof of
Theorem~\ref{th:G-parking=trees}.  
In other words, one would like to present a bijection between 
$G$-parking functions and oriented spanning trees of $G$.
There are several known bijections 
between the usual parking functions and trees.  One such bijection is 
relatively easy to construct.  There is a more elaborate bijection 
that maps parking functions $b$ with $b_1+\cdots+b_n=k$ to trees with 
$\binom{n}{2}-k$ inversions, see~\cite{Krew}.

\section{Square-free algebra}
\label{sec:CG}

Let $G$ be a graph on the set of vertices $0,\dots,n$.  
We will say that a subgraph $H\subset G$ of the graph $G$ 
is {\it slim\/} if the complement subgraph $G\setminus H$ is connected.
Let us associate commutative variables $\phi_e$, $e\in G$, with edges of the 
graph $G$, and let $\Phi_G$ be the algebra 
over $\K$ generated by the $\phi_e$ with the defining relations:
$$
\begin{array}{c}
(\phi_e)^2 = 0,\quad\textrm{for any edge $e$; } \\[.1in]
\ds\prod_{e\in H} \phi_e  = 0,\quad\textrm{for any non-slim
subgraph $H\subset G$.}
\end{array}
$$
Clearly, the square-free monomials $\phi_H=\prod_{e\in H} \phi_e$,
where $H$ ranges over all slim subgraphs in $G$,
form a linear basis of the algebra $\Phi_G$.
Thus the dimension of $\Phi_G$ is equal to the number of connected
subgraphs in $G$.

For $i=1,\dots,n$, let 
$$
X_i = \sum_{e\in G} c_{i,e}\, \phi_e,
$$
where 
$$
c_{i,e} = \left\{\begin{array}{rl} 
1 & \textrm{if } e = (i,j), \ i<j;\\
-1 & \textrm{if } e = (i,j), \ i>j;\\
0 & \textrm{otherwise.}
\end{array}\right.
$$
Define $\C_G$ as the subalgebra in $\Phi_G$ generated 
by the elements $X_1,\dots,X_n$.

Fix a linear ordering of edges of the graph $G$.
Recall that $N_G^k$ denotes the number of spanning trees of $G$
with external activity $k$, see Section~\ref{sec:power}.

\begin{theorem}
{\rm (1)} \
The dimension of the algebra $\C_G$ as a linear space over $\K$
equals the number of spanning trees in the graph $G$.

\medskip
\noindent
{\rm (2)} \ 
The dimension of the $k$-th graded component $\C_G^k$ 
of the algebra $\C_G$ equals the number 
$N_G^{|G|-n-k}$ of spanning trees of $G$ with external activity 
$|G| - n - k$.
\label{th:CG}
\end{theorem}

Recall that, for a nonempty subset $I\subset\{1,\dots,n\}$, 
$D_I=\sum_{i\in I,\,j\not\in I} a_{ij}$ is the number of edges
in $G$ that connect a vertex inside $I$ with a vertex outside of $I$,
see Section~\ref{sec:power}.

\begin{lemma}  For any nonempty subset $I\subset\{1,\dots,n\}$, 
the following relation holds in the algebra $\C_G$:
$$
\(\sum_{i\in I} X_i\)^{D_I} = 0.
$$
\label{lem:C<=B}
\end{lemma}

This lemma shows that the algebra $\C_G$ is a quotient of the algebra
$\B_G$.   We will eventually see that $\B_G=\C_G$, but let us pretend
that we do not know this yet.

\begin{proof}
Let $H_I\subset G$ be the subgraph that consists of all edges
that connect a vertex in $I$ with a vertex outside of $I$.
We have $\sum_{i\in I} X_i = \sum_{e\in H_I} \pm\phi_e$.
Thus $\(\sum_{i\in I} X_i\)^{D_I}=\pm\prod_{e\in H_I} \phi_e = 0$,
because $H_I$ is not a slim subgraph of $G$.
\end{proof}

Let $\S_G$ be the subspace in $\K[y_1,\dots,y_n]$ spanned by the elements
$$
\alpha_H = \prod_{e\in H} \alpha_e,
$$
for all slim subgraphs $H\subset G$,
where $\alpha_e = y_i - y_j$, for an edge $e=(i,j)$, $i<j$. 
Let $\S_G^k$ denote the $k$-th graded component of the space $\S_G$.

\begin{lemma}
For any graph $G$ and any $k$, we have $\dim \C_G^k = \dim \S_G^k$.
\label{lem:Gk}
\end{lemma}

\begin{proof}
Let $b_{H,a}$ be the coefficient of $\prod_{e\in H} \phi_e$ 
in the expansion of $X_{1}^{a_1}\cdots X_{n}^{a_n}$,
where $a=(a_1,\dots,a_n)$.
Then $\dim \C_G^k$ is equal to the rank of the matrix
$B=(b_{H,a})$, where $H$ ranges over all slim subgraphs in $G$
with $k$ edges and $a$ ranges over non-negative integer $n$-element sequences
with $a_1+\cdots+a_n=k$.  On the other hand,
$b_{H,a}$ is also equal to the coefficient of 
$y_1^{a_1}\cdots y_{n}^{a_n}$ in the expansion of $\alpha_H$.
Thus $\dim \S_G^k$ equals the rank of the same matrix $B=(b_{H,a})$.
\end{proof}

For a spanning tree $T$ in $G$, let $T^+$ denote the graph obtained
from $T$ by adding all externally active edges.  In virtue of 
Lemma~\ref{lem:Gk}, the following claim implies 
Theorem~\ref{th:CG}.

\begin{proposition}  The collection of elements 
$\alpha_{G\setminus T^+}$,
where $T$ ranges over all spanning trees of $G$, forms
a linear basis of the space $\S_G$.
\label{prop:SG}
\end{proposition}

Let us first prove a weaker version of this claim.
\begin{lemma}
The elements $\alpha_{G\setminus T^+}$,
where $T$ ranges over all spanning trees of $G$, spans
the space $\S_G$.
\label{lem:ST+}
\end{lemma}

\begin{proof}
Suppose not. 
Let $H$ be the lexicographically maximal slim subgraph of $G$ such 
that $\alpha_H$ cannot be expressed as a linear combination of
the $\alpha_{G\setminus T^+}$.
Then there exists a cycle $C=\{e_1,\dots,e_l\}\subset G$ with the minimal 
element $e_1$ such that $H\cap C=\{e_1\}$.
Then $\alpha_{e_1}$ is a linear combination of $\alpha_{e_2},
\dots,\alpha_{e_l}$.
Let $H_i$ be the graph obtained from $H$ by replacing the 
edge $e_1$ with $e_i$.
For $i=2,\dots,l$, the graph $H_i$ is a slim subgraph of $G$, which is 
lexicographically greater than $H_1$.
Then $\alpha_H$ can be expressed as a linear combination of
$\alpha_{H_2},\dots, \alpha_{H_l}$.  Contradiction.
\end{proof}

\begin{proof}[Proof of Proposition~\ref{prop:SG}]
Recall that $N_G$ denote the number of spanning trees in the graph $G$.
In view of Lemma~\ref{lem:ST+} it is enough to show that 
$\dim \S_G = N_G$.
We will prove this statement by induction on the number of edges in $G$.

If $G$ is a disconnected graph then it has no slim subgraphs
and $\dim \S_G = N_G=0$.
If $G$ is a tree then $\dim \S_G=N_G=1$. 
This establishes the base of induction.

Suppose that $G$ is a graph with at least one edge.  Pick an edge $e$ of $G$.
Let $G\backslash e$ be the graph obtained 
from $G$ by removing the edge $e$, also let $G/e$ be the graph obtained 
from $G$ by contracting the edge $e$.  
Then $N_G= N_{G\backslash e}+N_{G/e}$.  Indeed, for a spanning tree $T$ 
in $G$, we have either $e\not\in T$ or $e\in T$.  The former trees 
are exactly the spanning trees of $G\backslash e$.   The later trees 
are in a bijective correspondence with spanning trees of $G/e$.
This correspondence is given by contracting the edge $e$.
Assume by induction that the statement is true for both graphs 
$G\backslash e$ and $G/e$.

Let $\S_G'\subset \S_G$ be the span of the $\alpha_{H'}$
with slim subgraphs $H'\subset G$ such that $e\in H'$ and let 
$\S_G''\subset S_G$ be the span of the $\alpha_{H''}$
with slim subgraphs $H''\subset G$ such that $e\not\in H''$.
Then the space $\S_G$ is spanned by 
$\S_G'$ and $\S_G''$.  Thus 
\begin{equation}
\dim \S_G = \dim \S_G' + \dim \S_G'' - \dim (\S_G'\cap \S_G'').
\label{eq:dimSG}
\end{equation}
We have $\S_G'=(y_i-y_j)\, \S_{G\backslash e}$,
where $e=(i,j)$.  Thus $\dim \S_G' = \dim \S_{G\backslash e}$.
Let $p:f(y_1,\dots,y_n)\mapsto f(y_1,\dots,y_n)\,\mathrm{mod}\,(y_i-y_j)$ 
be the natural projection.
Then $p(\S_G'') = \S_{G/e}$ and $\S_G'\cap \S_G'' \subset \Ker(p)$.
Thus 
\begin{equation}
\dim \S_G'' = \dim\S_{G/e} + \dim\Ker(p) 
\geq \dim \S_{G/e} + \dim (\S_G'\cap \S_G'').
\label{eq:dimSG''}
\end{equation}
Combining~(\ref{eq:dimSG}) and~(\ref{eq:dimSG''}), we get
$$
\dim \S_G \geq \S_{G\backslash e} + \dim \S_{G/e}.
$$
By the induction hypothesis, the right-hand side of this expression
equals $N_{G\backslash e} + N_{G/e} = N_{G}$.
Thus $\dim \S_G \geq N_{G}$. 
On the other hand, Lemma~\ref{lem:ST+} implies that 
$\dim \S_G\leq N_{G}$.
Thus $\dim \S_G = N_{G}$, as needed.
\end{proof}

\section {Proof of Theorems~\ref{th:ABK}, \ref{th:dim-graded}, 
and~\ref{th:span}}
\label{sec:monote-proof}

Let $\{m_I\mid I\in \Sigma\}$ be a monotone monomial family,
and let $\I=\<m_I\mid I\in\Sigma\>$ be the corresponding monomial ideal
in $\K[x_1,\dots,x_n]$.

For a subset $I=\{i_1,\dots,i_r\}\in\{1,\dots,n\}$,
let $\K[x_I]= \K[x_{i_1},\dots,x_{i_k}]$, and
let $\MF_I$ denote the set of all monomials in the variables 
$x_i$, $i\in I$.  Also let $\MF=\MF_{\{1,\dots,n\}}$.
For $I\in \Sigma$, let $M_I=m_I\cdot \MF$ 
be the set of all monomials in $\K[x_1,\dots,x_n]$ divisible by $m_I$.
The standard monomial basis $R$ of the algebra 
$\A=\K[x_1,\dots,x_n]/\I$ is
the set of monomials 
$$
R=\MF\setminus\bigcup_{I\in\Sigma}M_I
$$ 
that survive in the algebra $\A$.

For $I, J\in\Sigma$, 
denote by $m_{J/I}$ the monomial obtained from 
$m_J$ by removing all $x_i$'s with $i\in I$, and let
$M_{J/I} = m_{J/I}\cdot \MF_{\overline{I}}$,
where $\overline{I}=\{1,\dots,n\}\setminus I$.
Let 
$\I_{I}$ be the monomial ideal in the polynomial ring 
$\K[x_{\overline{I}}]$ generated
by the monomials $\{m_{J/I}\mid J\in\Sigma, J\not\subset I\}$.
It follows from the monotonicity condition~(MM2) that
the ideal $\I_{I}$ is also generated by the set of monomials 
$\{m_{J/I}\mid J\in\Sigma, J\supsetneq I\}$.
Notice that $\I_I$ is also a monotone monomial ideal.
Let $R_I$ be the standard monomial basis of 
the algebra $\A_I=\K[x_{\overline{I}}]/\I_I$:
$$
R_I=\MF_{\overline{I}}
\setminus
\bigcup_{J\supsetneq I} M_{J/I}.
$$

\begin{proposition}
The polynomial ring $\K[x_1,\dots,x_n]$ decomposes into
the direct sum of subspaces:
$$
\K[x_1,\dots,x_n]=\<R\> \oplus
\bigoplus_{I\in\Sigma}
m_I\, \K[x_I]\,\<R_I\>,
$$
where $\<R\>$ and $\<R_I\>$ denote the linear
spans of monomials in $R$ and $R_I$, respectively.
\label{prop:decomp-m}
\end{proposition}

\begin{lemma} 
For any monomial $x^a=x_1^{a_1}\cdots x_n^{a_n}$ 
in the ideal $\I$ there is a unique
maximal by inclusion subset $J\in\Sigma$ such that $x^a\in M_I$.
\label{le:maximal} 
\end{lemma}

\begin{proof}
Let $\Sigma^\prime=\{I\in \Sigma\mid x^a\in M_I\}$.
If $I,J\in\Sigma^\prime$, then, according to condition (MM3),
there is an upper bound $K\in\Sigma$ of $I$ and $J$ such that
$M_I\cap M_J\subseteq M_K$.  Thus $M_K\in\Sigma'$.
This implies that $\Sigma'$ has a unique maximal element.
\end{proof}

\begin{proof}[Proof of Proposition~\ref{prop:decomp-m}]
For $I\in \Sigma$, let $\MM_I$ be the following set of monomials:
$$
\MM_I=M_I\setminus\bigcup_{J\supsetneq I}M_J,
$$
i.e., $\MM_I$ is the set of monomials 
$x^a\in \MF$ such that 
$I$ is the maximal by inclusion subset $I\in\Sigma$ with $x^a\in M_I$,
see Lemma~\ref{le:maximal}.
Thus the set of all monomials in $\K[x_1,\dots,x_n]$ 
decomposes into the disjoint union
\begin{equation}
\MF=R \ {\dot\cup} \  \dot{\bigcup_{I\in \Sigma}}\, \MM_I.
\label{eq:A-disjoint}
\end{equation}

Using monotonicity condition (MM2), we obtain, for $I\subsetneq J$, 
$$
M_I\setminus M_J = m_I\times \MF_I \times 
(\MF_{\overline{I}}\setminus M_{J/I}),
$$
where the notation ``$\times$'' means that every monomial in 
the left-hand side decomposes uniquely into the product of monomials.
Thus we have
\begin{equation}
\MM_I = \bigcap_{J\supsetneq I} (M_I\setminus M_J)
= m_I\times \MF_I\times 
\bigcap_{J\supsetneq I} (\MF_{\overline{I}}\setminus M_{J/I})
=m_I\times \MF_I\times R_I.
\label{eq:DI=mAR}
\end{equation}
Formulas~(\ref{eq:A-disjoint}) and~(\ref{eq:DI=mAR})
imply the required statement.
\end{proof}

Let $p_I$, $I\in\Sigma$, be a collection of polynomials
such that $p_I$ is an $I$-deformation of the monomial $m_I$. 
Remarkably, a similar statement  is valid for the polynomials $p_I$.

\begin{proposition}
The polynomial ring $\K[x_1,\dots,x_n]$ decomposes into
the direct sum of subspaces:
$$
\K[x_1,\dots,x_n]=\<R\> \oplus
\bigoplus _{I\in\Sigma}
p_I\, \K[x_I]\,\<R_I\> .
$$
\label{prop:decomp-p}
\end{proposition}

Proposition~\ref{prop:decomp-p}
immediately implies Theorem~\ref{th:span},
which says that the monomials in $R$ linearly span
the algebra $\B=\K[x_1,\dots, x_n]/\<p_I\mid I\in\Sigma\>$.

\begin {lemma}
Suppose that a polynomial $p\in \K[x_I]$ 
is an $I$-deformation of a monomial $m\in \K[x_I]$, 
see~{\rm(\ref{eq:genericity})}. 
Then for any polynomial 
$f\in\K[x_I]$ there exists a unique polynomial 
$g\in \K[x_I]$ such that the difference
$m\cdot f - p\cdot g$ contains no monomials divisible by $m$.
The map $f\mapsto g$ is one-to-one.
\label{lem:f-tilde-f}
\end{lemma}

\begin{proof} 
According to the generosity condition~(\ref{eq:genericity})
the polynomial $m\cdot f$, as well as any other polynomial in $\K[x_I]$, 
can be written uniquely in the form $m\cdot f = p\cdot g + r$,
where $g\in\K[x_I]$ and $r$ is in the linear span 
$\<R_m\>$ of monomials in $\K[x_I]$ not divisible by $m$.
This defines the map $f\mapsto g$.

On the other hand, for any $g\in\K[x_I]$ there exist unique
$f\in\K[x_I]$ and $r\in\<R_m\>$ such that $p\cdot g = m\cdot f -r$.
Thus the map $f\mapsto g$ is invertible.
The statement of the lemma follows.
\end{proof}

\begin{proof}[Proof of Proposition~\ref{prop:decomp-p}]
Pick any linear ordering $I_1,\dots,I_N$ of the set $\Sigma$
compatible with the inclusion relation, i.e, 
the inclusion $I_s\subset I_t$ implies that $s\leq t$.
Let $\Sigma_{(s)} = \{I_1,\dots,I_s\}$ and 
$\Sigma^{(s)} = \{I_s,\dots,I_N\}$ be initial and terminal 
intervals of $\Sigma$.

We will prove by induction on $N-s$ that 
the polynomial ring $\K[x_1,\dots,x_n]$ decomposes 
into the direct sum of subspaces
\begin{equation}
\K[x_1,\dots,x_n]=\<R\> \oplus
\bigoplus_{I'\in\Sigma_{(s)}} m_{I'}\,\K[x_{I'}]\,\<R_{I'}\> \oplus
\bigoplus_{I''\in\Sigma^{(s+1)}} p_{I''}\,\K[x_{I''}]\,\<R_{I''}\>. 
\label{eq:R-m-p}
\end{equation}
If $s=N$ then~(\ref{eq:R-m-p}) is true according to 
Proposition~\ref{prop:decomp-m}.  This gives the base of induction.

Assume by the induction hypothesis that~(\ref{eq:R-m-p}) holds
for some $s$ and derive the same statement for $s-1$.
Let $I=I_s$.
For a polynomial $\phi\in \K[x_1,\dots,x_n]$, write its unique
presentation
\begin{equation}
\phi=r +\sum_{I'\in \Sigma_{s-1}}
m_{I'}\cdot f_{I'}\cdot r_{I'}+
m_{I}\cdot f_{I}\cdot r_{I}+
\sum_{I''\in\Sigma^{s+1}} p_{I''}\cdot  f_{I''}\cdot r_{I''},
\label{eq:phi}
\end{equation}
where $r\in\< R\>$ and 
$f_{J}\in\K[x_J]$ and $r_{J}\in R_J$, for any $J\in\Sigma$.

Let $\tilde f_I\in\K[x_I]$ be the unique polynomial, 
provided by Lemma~\ref{lem:f-tilde-f}, such that the difference 
$d=m_I\cdot f_I-p_I \cdot \tilde f\in\K[x_I]$ contains no 
monomials divisible by $m_I$. 
Let $\psi\in\K[x_1,\dots,x_n]$ be the polynomial obtained from
$\phi$ by keeping all terms in~(\ref{eq:phi}) except 
for $m_{I}\cdot f_{I}\cdot r_{I}$
which we substitute  by the term $p_{I}\cdot \tilde f_{I}\cdot r_{I}$.
Then $\phi-\psi = d\cdot r_I$.
Pick any monomial $e$ in $d$.  
Remind that, according to~(\ref{eq:DI=mAR}),
$\MM_J$ is the set of all monomials in $m_J\,\K[x_J]\,\<R_J\>$.
If $e\cdot r_I\not\in \MM_J$ for all $J\in\Sigma$, 
then $e\cdot r_I\in\<R\>$.   
Otherwise, suppose that $e\cdot r_I\in \MM_J$ for some $J$.
If $J\not\subset I$, then $e\cdot r_I\in \MM_J\subset M_J$ implies that 
$r_I\in M_{J/I}$, which is impossible.
Thus $J\subset I$. Also $J\ne I$ because $e$ is not divisible by $m_I$.
This shows that
$$
\phi-\psi\in
\<R\>\oplus\bigoplus_{J\subsetneq I}\<\MM_J\> .
$$
Therefore, $\phi$ can be written as
\begin{equation}
\phi=\tilde r +\sum_{I'\in \Sigma_{s-1}}
m_{I'}\cdot \tilde f_{I'}\cdot \tilde r_{I'}+
p_{I}\cdot \tilde f_{I}\cdot r_{I}+
\sum_{I''\in\Sigma^{s+1}} p_{I''}\cdot f_{I''}\cdot r_{I''},
\label{eq:phi-next}
\end{equation}
where $\tilde r\in\< R\>$, 
$\tilde f_{I'}\in\K[x_{I'}]$, $\tilde r_{I'}\in R_{I'}$,
and $f_{I''}$ and $r_{I''}$ are the same as before.

Notice that all steps in the transformation of the 
presentation~(\ref{eq:phi}) to the presentation~(\ref{eq:phi-next}) 
are invertible.  Also if $p_I\cdot \tilde f_I\cdot r_I =0$ then
all summands in~(\ref{eq:phi}) and~(\ref{eq:phi-next}) coincide.
So, if at least the one of the summands in the presentation~(\ref{eq:phi-next}) 
of $\phi=0$ is nonzero, then we can also find a nonzero presentation
of the form~(\ref{eq:phi}) for $\phi=0$, which is impossible by the
induction hypothesis.  This shows that the presentation~(\ref{eq:phi-next}) 
of $\phi$ is unique.

This proves~(\ref{eq:R-m-p}).  For $s=0$ we obtain the claim
of Proposition~\ref{prop:decomp-p}.
\end{proof}

Finally we can put everything together and prove 
Theorems~\ref{th:ABK} and~\ref{th:dim-graded}.

\begin{proof}[Proof of Theorems~\ref{th:ABK} and~\ref{th:dim-graded}]
For a graph $G$, let $\A_G$, $\B_G$, 
and $\C_G$ be the algebras defined
in Sections~\ref{sec:G-parking}, \ref{sec:power}, and~\ref{sec:CG}.
Then we have the following termwise inequalities of Hilbert series
\begin{equation}
\Hilb\, \A_G \geq \Hilb\, \B_G \geq \Hilb\, \C_G.
\label{eq:A<B<C}
\end{equation}
The first inequality follows from Theorem~\ref{th:span}
because $\I_G$ is a monotone monomial ideal and, 
by Lemma~\ref{le:generic}, $\J_G$ is its deformation. 
The second inequality follows from Lemma~\ref{lem:C<=B}
that says that $\C_G$ is a quotient of $\B_G$.
Theorem~\ref{th:G-parking=trees} claims that $\dim \A_G=N_G$
is the number of spanning trees of the graph $G$.
On the other hand, by Theorem~\ref{th:CG}, $\dim \C_G=N_G$.
Thus all inequalities in~(\ref{eq:A<B<C}) are actually 
equalities.  Moreover, by Theorem~\ref{th:CG}, the dimensions
of $k$-th graded components are equal to 
$$
\dim \A_G^k = \dim \B_G^k = \dim \C_G^k = N_G^{|G|-n-k},
$$
the number of spanning trees of $G$ with external activity 
$|G|-n-k$.
\end{proof}

\begin{corollary}
The algebras $\B_G$ and $\C_G$ are isomorphic.
\end{corollary}

\section{Algebras related to forests}
\label{sec:forests}

Definitions of the algebras $\B_G$ and $\C_G$ and the proof
of Theorem~\ref{th:CG} are similar to
constructions from~\cite{PSS}.
Let us briefly review some results from~\cite{PSS}.

Let $G$ be a graph on the vertices $0,\dots,n$.  Let $\hat\J_G$
be the ideal in $\K[x_1,\dots,x_n]$ generated by the polynomials
$$
\hat p_I = \left(\sum_{i\in I} x_i\right)^{D_I+1},
$$
where $I$ ranges over all nonempty subsets in $\{1,\dots,n\}$ and
the number $D_I$ is the same as in Section~\ref{sec:power},
cf.~(\ref{eq:pI}).
Let $\hat \B_G= \K[x_1,\dots,x_n]/\hat\J_G$.

Let $\hat \Phi_G$ be the commutative algebra generated by 
the variables $\hat \phi_e$, $e\in G$, with the 
defining relations:
$$
(\hat\phi_e)^2=0, \qquad\textrm{for any edge } e.
$$
And let $\hat\C_G$ be the subalgebra of $\hat\Phi_G$ generated
by the elements
$$
\hat X_i = \sum_{e\in G} e_{i,e} \hat\phi_e,
$$
for $i=1,\dots,n$, cf.\ Section~\ref{sec:CG}.

A {\it forest\/} is a graph without cycles. 
The connected components of a forest are trees.  A {\it subforest\/} 
in a graph $G$ is a subgraph $F\subset G$ without cycles.
Fix a linear order of edges of $G$.  An edge $e\in G\setminus F$ is called 
{\it externally active\/} for a forest $F$ 
if there exists a cycle $C$ in $G$ such that
$e$ is the minimal element of $C$ and $(C\setminus\{e\})\subset F$. 
The {\it external activity\/} of $F$ is the number of externally active
edges for $F$.

\begin{theorem} {\rm \cite{PSS}} \ 
The algebras $\hat \B_G$ and $\hat \C_G$ are isomorphic to each other.
Their dimension is equal to the number of subforests in the graph $G$.

The dimension $\dim \hat\B_G^k$ of the $k$-th graded component 
of the algebra $\hat\B_G$ equals the number of subforests $F$ of $G$
with external activity $|G|-|F|-k$.
\end{theorem}

In~\cite{PSS0} we investigated the algebra
$\hat\B_G$ for the graph $G=K_{n+1}$.
Let $\hat\I_n=\<\hat m_I\>$ and 
$\hat\J_n=\<\hat p_I\>$ be two ideals in the polynomial ring 
$\K[x_1,\dots,x_n]$ generated by the monomials $\hat m_I$ 
and the polynomials $\hat p_I$, correspondingly, given by
$$
\begin{array}{c}
\hat m_I=(x_{i_1}\cdots x_{i_r})^{n-r+1} x_{i_1},\\[.1in]
\hat p_I=(x_{i_1}+\cdots+x_{i_r})^{r(n-r+1)+1},
\end{array}
$$
where $I=\{i_1<\dots<i_r\}$,
ranges over nonempty subsets of $\{1,\dots,n\}$,
cf.~Subsection~\ref{ssec:(n+1)^{n-1}}.
Notice that $\hat\I_n$ is a monotone monomial ideal
and $\hat\J_n$ is its deformation.
Let $\hat\A_n = \K[x_1,\dots,x_n]/\hat\I_n$ and 
$\hat\B_n = \K[x_1,\dots,x_n]/\hat\J_n$.

Let us say that a non-negative integer sequence 
$b=(b_1,\dots,b_n)$ is an {\it almost parking function} of size $n$ 
if the monomial 
$x^b =x_1^{b_1}\cdots x_n^{b_n}$ does not belong to the ideal~$\hat\I_n$.
Clearly the class of almost parking functions includes usual parking
functions.

For a forest $F$ on the vertices $0,\dots,n$,
an {\it inversion\/} is a pair of vertices labelled $i$ and $j$
such that $i>j$ and the vertex $i$ belong to the path
in $F$ that joins the vertex $j$ with the minimal vertex in its connected
component.

\begin{theorem}
{\rm \cite{PSS, PSS0}} \
The algebras $\hat\A_n$ and $\hat\B_n$ have the same Hilbert series.
The dimension of these algebras is equal to the number of forests
on $n+1$ vertices.

Moreover, the dimension $\dim \hat \A_n^k = \dim \hat\B_n^k$ of the 
$k$-th graded components of the algebras $\hat\A_n$ and $\hat\B_n$ 
is equal to
\begin{enumerate}
\item[(A)]
the number of almost parking functions $b$ of size $n$ such that
$\sum_{i=1}^n b_i = k$;
\item[(B)]
the number of forests on $n+1$ vertices with external activity 
$\binom{n}{2} -k$;
\item[(C)]
the number of forests on $n+1$ vertices with 
$\binom{n}{2} -k$ inversions.
\end{enumerate}

The images of the monomials $x^b$, where $b$ ranges over almost
parking functions of size $n$, form linear bases in both algebras
$\hat \A_n$ and $\hat\B_n$.
\label{th:A=B-hat}
\end{theorem}

Theorem~\ref{th:A=B-hat}, first stated in~\cite{PSS0},
follows from results of~\cite{PSS}.  
The algebra $\hat\B_n$ is the algebra generated by curvature forms 
on the complete flag manifold.  It was introduced in an attempt to 
lift Schubert calculus on the level of differential forms, 
see~\cite{PSS,PSS0,SS}.  This example related to Schubert calculus 
was our original motivation.

\section{$\rho$-algebras and $\rho$-parking functions}
\label{sec:rho-algebras}

We conclude the paper with a discussion of a special class 
of monotone monomial ideals and their deformations.

Let $\rho=(\rho_1,\dots,\rho_n)$
be a weakly decreasing sequence of non-negative integers, 
called a {\it degree function.} 
Let $\I_{\rho}=\<m_I\>$ and $\J_\rho=\<p_I\>$
be the ideals the ring $\K[x_1,\dots,x_n]$ 
generated by the monomials $m_I$ and the polynomials $p_I$, 
correspondingly, given by 
$$
\begin{array}{c}
\ds m_I=(x_{i_1}\cdots x_{i_r})^{\rho_r},\\[.1in]
\ds p_I=(x_{i_1}+\cdots+x_{i_r})^{r\cdot \rho_r},
\end{array}
$$
where in both cases $I=\{i_1,\dots,i_r\}$ runs over all
nonempty subsets of $\{1,\dots,n\}$.
Let $\A_{\rho} = \K[x_1,\dots,x_n]/\I_{\rho}$ and 
$\B_{\rho} = \K[x_1,\dots,x_n]/\J_{\rho}$.

Let us say that a non-negative integer sequence $b=(b_1,\dots,b_n)$
is a {\it $\rho$-parking function\/} if the monomial
$x_1^{b_1}\cdots x_n^{b_n}$ does not belong to the ideal $\I_\rho$.
More explicitly, this condition can be reformulated as follows.
A non-negative integer sequence $b=(b_1,\dots,b_n)$ 
is a $\rho$-parking function if and only if,
for $r=1,\dots,n$, we have
$$
\#\{i\mid b_i< \rho_{n-r+1}\} \geq r.
$$
This condition can also be formulated in terms of 
the increasing rearrangement $c_1\leq \cdots\leq c_n$ of the elements 
of $b$ as $c_i<\rho_{n+1-i}$.
The $\rho$-parking functions were studied in~\cite{PiSt} and in~\cite{Yan}.
They also appeared under a different name in~\cite{PP}. 
Notice that $(n,\dots,1)$-parking functions are exactly
the usual parking functions of size $n$.

The monomials $x^b$, where $b$ ranges over $\rho$-parking
functions, form the standard monomial basis of the algebra $\A_\rho$.
Thus the Hilbert series of the algebra $\A_\rho$ equals
$$
\Hilb\, \A_\rho = \sum_{b} q^{b_1+\cdots+b_n},
$$ 
where the sum is over $\rho$-parking functions. 
The dimension $\dim\A_\rho$ of this algebra 
is equal to the number of $\rho$-parking functions.

Theorem~\ref{th:span} specializes to the following statement. 

\begin{corollary} 
The monomials $x^b$, where $b$ ranges over $\rho$-parking
functions, linearly span the algebra $\B_\rho$.
Thus we have the termwise inequality of Hilbert series:
$$
\Hilb\,\A_\rho\ge \Hilb\,\B_\rho.
$$ 
\end{corollary}

It would be interesting to describe the class of degree functions $\rho$
such that $\Hilb\,\A_\rho = \Hilb\,\B_\rho$.
If $\rho_r= l+k(n-r)$ is a linear degree function then,
according to Corollary~\ref{cor:A=Bnkl},
the Hilbert series of $\A_\rho$ and $\B_\rho$ are equal to each other
and
$$
\dim \A_\rho = \dim \B_\rho = l\,(l+kn)^{n-1}. 
$$
For $n=3$, Schenck~\cite{Schenck} gave another proof of 
this fact using ideals of fatpoints.

Let us say that a degree function $\rho$ is {\it almost linear\/} if there
exists an integer $k$ such that $\rho_i-\rho_{i+1}$ equals either $k$ or
$k+1$, for $i=1,\ldots, n-1$.
Computer experiments show that 
the equality $\Hilb\,\A_{\rho}= \Hilb\,\B_{\rho}$
often holds for almost linear degree functions $\rho$.
The table below  lists some almost linear degree functions,
for which the equality $\Hilb\,\A_\rho  = \Hilb\,\B_\rho$ holds.

\begin{table}[ht]
\footnotesize
\begin{tabular}{|@{} c  @{}  | c @{} c ||@{} c @{} c @{}|   c |@{}}
\hline
\quad \quad $\rho$ \quad \quad & &
$\dim \A_\rho$ & \quad  \quad $\rho$ \quad\quad & & $\dim \A_\rho$ \\
\hline
$(4,2,1)$ && 25  & $(8,5,3)$  && 306 \\
$(8,5,1)$ && 142 & $(8,6,3)$  && 351 \\
$(6,4,3)$ && 153 & $(11,7,2)$ && 506 \\
$(9,5,2)$ && 290 & $(12,8,3)$ && 855 \\
\hline
$(6,4,3,2)$ && 632  & $(8,7,5,3)$  && 3021  \\
$(9,6,4,2)$ && 2512 & $(9,7,6,5)$  && 4925  \\
$(8,6,5,3)$ && 2643 & $(11,8,6,3)$ && 7587  \\
$(8,6,5,4)$ && 2832 & $(12,9,7,4)$ && 12460 \\
\hline
\ $(9,8,6,4,2)$ \ && 31472 & \ $(10,9,7,5,3)$ \ && 65718 \\
\hline
\end{tabular}
\normalsize
\end{table}

On the other hand, the equality of Hilbert series fails for 
the almost linear degree functions $\rho= (9,6,3,1)$ and $\rho=(9,7,5,4,3)$. 
We do not know an example when $\Hilb\,\A_{\rho}= \Hilb\,\B_{\rho}$ and
$\rho$ is not almost linear.

The ideal $\I_\rho$ is a {\it strictly\/} monotone monomial ideal
provided that the degree function is strictly decreasing 
$\rho_1>\cdots > \rho_n>0$.
Corollary~\ref{cor:resolution2} gives a minimal free resolution for
this ideal.  Recall that $S=\K[x_1,\dots,x_n]$ and $S(-d)$ is the
free $\Z$-graded $S$-module of rank 1 with generator of degree $d$.

\begin{corollary} Let $\rho$ be a degree function such that
$\rho_1>\cdots > \rho_n>0$.
The ideal $\I_\rho$ has a minimal free resolution of the form
$$
\cdots 
\longrightarrow 
C_3\longrightarrow 
C_2\longrightarrow 
C_1\longrightarrow 
C_0=S\longrightarrow 
S/\I_\rho\longrightarrow 0,
$$
with
$$
C_k = \bigoplus_{l_1,\dots,l_k} S(-d(l_1,\dots,l_k))^{\binom{n}{l_1,\cdots,l_k}}\,,
$$
where the direct sum is over 
$l_1,\dots,l_k\geq 1$ such that $l_1+\cdots +l_k\leq n$,
$$
d(l_1,\dots,l_k) = l_1\,\rho(l_1)+l_2\,\rho(l_1+l_2)+\cdots
+l_k\,\rho(l_1+\cdots+l_k),
$$
and $\binom{n}{l_1,\cdots,l_k}=
\frac{n!}{l_1!\cdots l_k! (n-l_1-\cdots - l_k)!}$
is the multinomial coefficient.
\label{cor:rho-resolution}
\end{corollary}

Conjecture~\ref{conj:Betti=Betti} imply that if 
$\dim \A_\rho=\dim \B_\rho$ then the ideals $\I_\rho$
and $\J_\rho$ have the same graded Betti numbers.
It is already a nontrivial open problem to prove (or disprove) 
that the graded Betti numbers of the ideal $\J_\rho$, 
for a linear degree function $\rho$, are given by
the expression in Corollary~\ref{cor:rho-resolution}.
Schenck's results~\cite{Schenck} for $n=3$ seem to support
this conjecture.

Proposition~\ref{prop:monotone-GF} specializes to an expression for 
the Hilbert series $\Hilb\,\A_\rho$ with alternating signs.
Actually, in this case it is possible to give a simpler subtraction-free 
expression for the Hilbert series.
The following statement is a slight enhancement of a result of 
Pitman and Stanley, who gave a formula for the number of $\rho$-parking
functions.

\begin{proposition}
{\rm cf.~\cite[Theorem~11]{PiSt}} \
The Hilbert series of the algebra $\A_\rho$ equals
\begin{equation}
\Hilb\,\A_\rho = \sum_{a} \,\prod_{i=1}^n 
\frac{q^{\rho_{n-a_i+1}} - q^{\rho_{n-a_i}}}{1-q},
\label{eq:HArho}
\end{equation}
where the sum is over $(n+1)^{n-1}$ usual parking functions
$a=(a_1,\dots,a_n)$ of size $n$.  Here we assume that $\rho_{n+1}=0$.
Thus the dimension of $\A_\rho$, which is the number of 
$\rho$-parking functions, is given by the following polynomial
in $\rho_1,\dots,\rho_n$:
$$
\dim\A_\rho = \sum_{a} \,\prod_{i=1}^n 
(\rho_{n-a_i} - \rho_{n-a_i+1}),
$$
where again the sum is over usual parking functions of size $n$.
\end{proposition}

\begin{proof}
For $i=0,\dots,n$, let $Z_i$ be the interval of integers
$Z_i=[\rho_{n-i+1},\rho_{n-i}[$, 
where we assume that $\rho_0=+\infty$ and $\rho_{n+1}=0$.
Then the set of positive integers is the disjoint union of 
$Z_0,\dots,Z_n$.
Let $f:b\mapsto a$ be the map that sends
a positive integer sequence $b=(b_1,\dots,b_n)$
to the sequence $a=(a_1,\dots,a_n)$ such that $b_i\in Z_{a_i}$ 
for $i=1,\dots,n$.
Then $b$ is a $\rho$-parking function if and only if $a$ is
a usual parking function of size $n$.
Fix a parking function $a$ of size $n$.  Then 
$$
\sum_{b:\,f(b)=a} q^{b_1+\cdots+b_n} =
\prod_{i=1}^n \,\sum_{b_i\in Z_{a_i}} q^{b_i}.
$$
is exactly the summand in~(\ref{eq:HArho}).
\end{proof}

For example, the Hilbert series of $\A_{\rho}$, for $n=2$ and $n=3$, 
are given by 
$$
\begin{array}{l}
\Hilb\,\A_{(\rho_1,\rho_2)}(q) = 
[\rho_2]^2 + 2q^{\rho_2}[\rho_1-\rho_2]\,[\rho_2]\,,\\[.1in]
\Hilb\,\A_{(\rho_1,\rho_2,\rho_3)}(q) = 
[\rho_3]^3 
+ 3q^{\rho_3}[\rho_2-\rho_3]\,[\rho_3]^2
+ 3q^{2\rho_3}[\rho_2-\rho_3]^2\,[\rho_3]+\\[.1in]
\qquad\qquad\qquad
+ 3q^{\rho_2}[\rho_1-\rho_2]\,[\rho_3]^2
+ 6q^{\rho_2+\rho_3}[\rho_1-\rho_2]\,[\rho_2-\rho_3]\,[\rho_3]\,,
\end{array}
$$
where $[s] = 1+q+\cdots+q^{s-1}$ denotes the $q$-analogue of an integer $s$.

Finally, we formulate a theorem that gives a combinatorial interpretation
of the value of the Hilbert series $\Hilb\,\A_\rho$ at $q=-1$.
This theorem follows from results of~\cite{PP} 
on $\rho$-parking functions.

\begin{theorem} {\rm \cite{PP}} \ 
The number $(-1)^{\rho_1+\cdots+\rho_n-n}\,\Hilb\,\A_\rho(-1)$ 
equals the number of permutations $\sigma_1,\dots,\sigma_n$ of 
$1,\dots,n$ such that
$$
\sigma_1 \vee^{\rho_1} \sigma_2 \vee^{\rho_2}\cdots 
\vee^{\rho_{n-1}} \sigma_n\vee^{\rho_n} 0\,,
$$
where the 
notation $a\vee^k b$ means that $a<b$ for even $k$ and $a>b$ for odd $k$.
In particular, $\Hilb\,\A_\rho(-1)$ is zero if and only if 
$\rho_n$ is even.
\end{theorem}

This theorem basically says that $\Hilb\,\A_\rho(-1)$
is either $0$ or plus/minus the number of permutations with 
prescribed descent positions.

In the case of usual parking functions of size $n$, 
i.e., for $\rho=(n,\dots,1)$, this theorem amounts to the well-known 
result of Kreweras~\cite{Krew} that the value of the inversion polynomial
$I_n(-1)=(-1)^{\binom{n}{2}}\, \Hilb\,\A_{(n,\dots,1)}(-1)$ 
is the number of alternating permutations of size $n$.

\section{Appendix: Abelian sandpile model}
\label{sec:sandpiles}

In this appendix we discuss the {\it abelian sandpile model}, also
known as the {\it chip-firing game}.  It was introduced by
Dhar~\cite{Dhar} and was studied by many authors.  We review the
sandpile model for a class of toppling matrices introduced by Gabrielov~\cite{Gab2},
which is more general than in~\cite{Dhar}. Then we show how $G$-parking
functions from Section~\ref{sec:G-parking} are related to this model.

Let $\Delta=(\Delta_{ij})_{1\leq i,j\leq n}$ be an integer 
$n\times n$-matrix.  We say that $\Delta$ is a {\it toppling matrix\/} 
if it satisfies the following two conditions:
\begin{equation}
\Delta_{ij}\leq 0 \textrm{, for $i\ne j$;}\quad\textrm{and}\quad
\text{there exists a vector $h>0$ such that $\Delta\,h>0$.}
\label{eq:toppling-conditions}
\end{equation}
Here the notation $h>0$ means that all coordinates of $h$ are strictly positive.
Notice that conditions~(\ref{eq:toppling-conditions}) imply that
$\Delta_{ii}>0$ for any $i$.
These matrices appeared in~\cite{Gab2} under the name {\it
avalanche-finite redistribution matrices}.

Let us list some properties of toppling matrices.
Recall that $L_G$ denotes the truncated Laplace matrix
that corresponds to a digraph $G$ on the vertices $0,\dots,n$,
see Equation~{\rm(\ref{eq:Laplace})} in Section~\ref{sec:G-parking}.

\begin{proposition} {\rm cf.~\cite{Gab2}} \ 
{\bf 1.} A matrix $\Delta$ is a toppling matrix if and only if its 
transposed matrix $\Delta^T$ is a toppling matrix.

\smallskip
\noindent
{\bf 2.} Every integer matrix $\Delta$ such that
\begin{equation}
\qquad \Delta_{ij}\leq 0 \textrm{ for $i\ne j$;}\qquad
\sum_j \Delta_{ij}\geq 0 \textrm{ for all $i$;} \qquad
\det\Delta\ne 0
\label{eq:toppling=laplace}
\end{equation}
is a toppling matrix.
Equivalently, the truncated Laplace matrix $\Delta=L_G$ corresponding
to a digraph $G$ with at least one oriented spanning tree
is a toppling matrix.

\smallskip
\noindent
{\bf 3.} If $\Delta$ is a toppling matrix then all principal minors of
$\Delta$ are strictly positive.

\smallskip
\noindent
{\bf 4.}  If $\Delta$ is a symmetric integer matrix
with non-positive off-diagonal entries, then $\Delta$ is a toppling matrix 
if and only if it is positive-definite.
\label{prop:toppling-properties}
\end{proposition}

\begin{proof}
{\bf 1.} 
This claim follows from~\cite[Theorem~1.5]{Gab2}.
It also follows from the result of~\cite[Theorem~4.3]{Kac}, 
obtained  for classification of generalized Cartan matrices.

%implies the following fact:  
%Let $B=(b_{ij})$ be a real matrix such that $b_{ij}<0$ for $i\ne j$.  If there
%exists a vector $h>0$ such that $B\,h>0$ then the same property holds for the
%transposed matrix $B^T$.
%
%Suppose that $\Delta$ is a toppling matrix and $h>0$ is a vector such
%that $\Delta\,h>0$. Let $J$ be the $n\times n$ matrix with all entries
%equal $1$.  Then we can find a sufficiently small $\epsilon>0$ 
%such that $\tilde \Delta = \Delta-\epsilon J$ satisfies
%$\tilde \Delta\,h>0$.  The matrix $\tilde \Delta$ satisfies Kac's conditions.  
%Thus there is a vector $h'>0$ such that
%$\tilde\Delta^T\, h'>0$.  Then $\Delta^T\,h'\geq \tilde\Delta^T\,h'>0$  
%and $\Delta^T$ is a toppling matrix.

\smallskip
{\bf 2.}
Conditions~(\ref{eq:toppling=laplace}) are equivalent to the statement 
that $\Delta=L_G$ is the truncated Laplace matrix for some digraph $G$ 
with at least one oriented spanning tree, see the Matrix-Tree Theorem,
Equation~(\ref{eq:matrix-tree}) in Section~\ref{sec:G-parking}.
Let $\mathrm{dist}(i)$ be the length of the shortest directed path
in the digraph $G$ from the vertex $i$ to the root $0$, and let
$h(\epsilon) = (h_1,\dots,h_n)^T$, 
where $h_i = 1 - \epsilon^{\mathrm{dist}(i)}$.
Then $\Delta \,h(\epsilon)>0$ for sufficiently small $\epsilon>0$.
Indeed, the $i$-th coordinate of the vector $\Delta \,h(\epsilon)$
is 
$$a_{i0} (1 - \epsilon) + 
\sum_{j\ne0,i} a_{ij} (\epsilon^{\mathrm{dist}(j)}
-\epsilon^{\mathrm{dist}(i)}),
$$
where the $a_{ij}$ are the entries the adjacency matrix of $G$.
The leading term of this expression 
has order of $\epsilon^{\mathrm{dist}(i)-1}$ and is strictly positive.

\smallskip
{\bf 3.} 
The fact that $\det\Delta>0$ is given in~\cite[Proposition~1.12]{Gab2}.
Let us show that it also easily follows from the Matrix-Tree theorem.
Let $\Delta$ be a toppling matrix and $h=(h_1,\dots,h_n)^T>0$ 
be an integer vector such that $\Delta\,h>0$.  Then all row sums of the
matrix $\tilde \Delta = \Delta\cdot \mathrm{diag}(h_1,\dots,h_n)$ are positive.
This means $\tilde\Delta=L_G$ is the truncated Laplace 
matrix for some digraph $G$.  The $i$-th row sum of $\tilde\Delta$ 
is the number of edges in $G$ connecting the vertex
$i$ with the root $0$.  According to the Matrix-Tree 
Theorem,
the determinant $\det \tilde\Delta$ is the number of oriented spanning 
trees in the digraph $G$.  This number is positive, because each vertex 
is connected with the root by an edge in $G$. Thus 
$\det \Delta = (h_1\cdots h_n)^{-1}\det \tilde \Delta>0$.
Any principal minor of $\tilde \Delta$ also has positive row sums.
The same argument holds for the minors.

\smallskip
{\bf 4.}  This claim follows from~\cite[Lemma~4.5]{Kac}.
\end{proof}

Let us now fix a toppling matrix $\Delta$, and let
$\Delta_i=(\Delta_{i1},\dots,\Delta_{in})$ be the $i$-th row of $\Delta$.
A {\it configuration\/} $u=(u_1,\dots,u_n)$ is a vector of non-negative 
integers.  In the sandpile model, the number $u_i$ is interpreted as the 
the number of particles, or grains of sand, at site $i=1,\dots,n$.  
A site $i$ is {\it critical\/} if $u_i\geq \Delta_{ii}$.
A {\it toppling\/} at a critical site $i$ consists in subtraction the 
vector $\Delta_i$ from the vector $u$.
In other words, toppling at site $i$ decreases $u_i$ by $\Delta_{ii}$
particles and increases $u_j$ by $-\Delta_{ij}$ particles, for all $j\ne i$.
A configuration $u$ is called {\it stable\/} if no 
toppling is possible, i.e., $0\leq u_i<\Delta_{ii}$ for all sites $i$.

Dhar~\cite{Dhar} assumed that the toppling matrix $\Delta$ has
non-negative row sums, i.e., he assumed that $\Delta=L_G$ satisfies
conditions~(\ref{eq:toppling=laplace}).  In this case a toppling cannot
increase the total number of particles.  When a toppling occurs, some
of the particles at site $i$ are distributed among the neighboring
sites and some particles are removed from the system.  While this
condition is important from the physical point of view, it is not
really necessary for the following algebraic constructions,
cf.~Gabrielov~\cite{Gab2}.  Moreover,
there are interesting examples, for which this condition fails.  We
still attribute the following results to Dhar even though we will not
assume that $\Delta$ has non-negative row sums.  The proofs that we
include for completeness sake are close to the proofs from~\cite{Dhar}.

\begin{lemma} {\rm \cite{Gab2}, cf.~\cite{Dhar}} \ 
Every configuration can be transformed into a stable configuration
by a sequence of topplings.  This stable configuration 
does not depend on the order in which topplings are performed.
\end{lemma}

\begin{proof}
Conditions~(\ref{eq:toppling-conditions}) for a toppling matrix imply 
that there exists a vector $h>0$ such that $(h,\Delta_i)>0$ for any $i$.  
For a configuration $u$, the value $h(u)$ is non-negative 
and every toppling strictly decreases this value.  
Thus, after at most $h(u)/\min(h,\Delta_i)$  topplings, 
the configuration $u$ transforms into a stable configuration.

If an unstable configuration $u$ has two critical sites $i$ and $j$,
then $j$ is still a critical site for $u-\Delta_i$. 
Thus is it possible to perform a toppling at site $i$ followed by a toppling 
at site $j$ producing the configuration $u-\Delta_i-\Delta_j$.
This operation is symmetric in $i$ and $j$. 
Using this argument repeatedly, we deduce that the final stable 
configuration does not depend on the order of topplings.
\end{proof}

The {\it avalanche operators\/} $A_1,\dots,A_n$ map the set
of stable configurations to itself.  The operator $A_i$ is given 
by adding 1 particle at site $i$, i.e., increasing $u_i$ by $1$, 
and then performing a sequence of topplings that lead 
to a new stable configuration.

\begin{lemma} {\rm \cite{Dhar}} \ 
The avalanche operators $A_1, \dots,A_n$ commute pairwise. 
\label{lem:avalanche-commute}
\end{lemma}

\begin{proof} The stable configuration $A_i A_j u$ is obtained from $u$
by adding a particle at site $j$ then performing a sequence of topplings then
adding a particle at site $i$ and performing another sequence of topplings.
If we first add two particles at sites $i$ and $j$, then all topplings
in these two sequences are still possible and lead 
to the same stable configuration.
This shows that $A_i$ and $A_j$ commute.
\end{proof}

The {\it abelian sandpile model\/} is the random walk on the set of stable 
configurations that is given by picking a site $i$ at random 
with some probability $p_i>0$ and performing the avalanche 
operator $A_i$.
Informally, we can describe it as the model where we drop a grain
of sand at a random site and allow the system to settle to a
stable configuration.

Dhar described the steady state of this random walk.
A stable configuration $u$ is called {\it recurrent\/} if there are 
positive integers $c_i$ such that $A_i^{c_i} u = u$ for all $i$.
Let $\R$ denote the set of recurrent configurations.
The commutativity of the avalanche operators $A_i$ 
implies that the set $\R$ is closed under the action of these operators.  
Moreover, the operators $A_i$ are invertible on the set $\R$.  
Indeed, $A_i^{-1} u$ can be defined as $a_i^{c_i-1} u$ 
for a recurrent configuration $u$.  According to the theory of Markov 
chains all recurrent configurations have the same nonzero probability 
of occurrence in the steady state and all non-recurrent configurations 
have zero probability.

The {\it sandpile group\/} $SG$, also known as the 
{\it critical group}, is the finite abelian group generated 
by the avalanche operators $A_1,\dots,A_n$ acting on the set $\R$.

\begin{theorem}  {\rm \cite{Dhar}} \  
The sandpile group is isomorphic to the quotient of the integer lattice
$SG\simeq \Z^n/\<\Delta\>$,
where $\<\Delta\>=\Z\Delta_1\oplus\cdots\oplus \Z\Delta_n$
is the sublattice in $\Z^n$ spanned by the vectors $\Delta_i$.
The order of this group is equal to the number of recurrent configurations 
and is given by $|SG|=|\R|=\det \Delta$.
\label{th:order-sandpiles}
\end{theorem}

\begin{proof}
Since there are finitely many recurrent configurations, we may assume
that the numbers $c_i$ are the same for all recurrent configurations.
For a recurrent configuration $u\in\R$ and an integer vector 
$v=(v_1,\dots,v_n)$, let $A^v u = A_1^{v_1}\cdots A_n^{v_n} u$.
Then $u = A^{u-v} v$ for any $u,v\in\R$.  Indeed, 
the configuration $A^{u-v} v$ is given by performing topplings to the 
configuration $(u-v + N c) + v = u + N c$, where $c=(c_1,\dots,c_n)$
and $N$ an integer large enough to make the vector $u-v+Nc$ positive.  
The result of these toppling equals $A^{Nc}u = u$.  
This shows that $SG$ acts transitively on $\R$.
If an element of $SG$ stabilizes a configuration $u\in\R$ then, by
transitivity, it stabilizes any other element of $\R$ and is the
identity in $SG$.  Thus the order of the sandpile group $SG$
equals to $|\R|$.  The bijection between $SG$ and $\R$ is given by
$A^v\mapsto A^v\cdot u_*$, where $u_*$ is any fixed element of $\R$.

If we add $\Delta_{ii}$ particles at site $i$
to a configuration $u$ and perform a toppling at the 
(unstable) site $i$, the result will 
be the same as adding $-\Delta_{ij}$ particles at
all other sites $j\ne i$.  Thus
$$
A_i^{\Delta_{ii}} = \prod_{j\ne i} A_j^{-\Delta_{ij}},
\textrm{ or equivalently, }
A^{\Delta_i} = 1\textrm{ for any } i.
$$
On the other hand, $A^v\ne 1$ if $v\not\in \<\Delta\>$, since topplings
are given by subtraction of the vectors $\Delta_i$
and $A^v u\in u + v + \<\Delta\>$.  
This shows that $A^v=1$ if and only if $v\in\<\Delta\>$
and the map $v\mapsto A^v$ is an isomorphism between
the sandpile group $SG$ and the quotient $\Z^n/\Z\Delta$.

Finally, the order of $\Z^n/\<\Delta\>$ equals $\det\Delta$.
\end{proof}

Dhar suggested a more explicit characterization of the set $\R$
of recurrent configurations.
Let us say that a configuration  $u$ is {\it allowed\/}
if for any nonempty subset $I$ of sites there exists $j\in I$ 
such that 
$$
u_j\geq \sum_{i\in I\setminus\{j\}} (-\Delta_{ij}).
$$ 

\begin{proposition} {\rm \cite{Dhar}} \  
Every recurrent configuration is allowed.
\label{prop:recurrent->allowed}
\end{proposition}

\begin{proof}
Let $u$ be a recurrent configuration.
Then $A^c\,u = u$, where $c=(c_1,\dots,c_n)$ and $c_i>0$ for all $i$.
This means that there exists a sequence of sits $i_1,\dots,i_k$ such that
(i) $\Delta_{i_1}+\cdots + \Delta_{i_k} = c$; and  
(ii) $u+\Delta_{i_1}+\cdots + \Delta_{i_r}\geq 0$ for any $r=1,\dots,k$.
Since all coordinates of $\Delta_{i}$, except the $i$-th coordinate, 
are non-positive and $c>0$, condition (i) implies that 
the sequence $i_1,\dots,i_k$ contains all sites $1,\dots,n$ at least once.

Let us say that a configuration $v$ is {\it $I$-forbidden\/}, for some
subset $I$ of sites, if 
$$
0\leq v_j < \sum_{i\in I\setminus\{j\}} (-\Delta_{ij}),
$$
for all $j\in I$.  If $v$ is $I$-forbidden and $v+\Delta_i\geq 0$, 
for some site $i$, then the configuration $v+\Delta_i$ is 
$I\setminus\{i\}$-forbidden.
Also notice that there are no $\emptyset$-forbidden
configurations. 

Suppose that the recurrent configuration $u$ is not allowed.
Then $u$ is $I$-forbidden of some subset $I$.
We obtain by induction on $r=0,\dots,k$ that
the configuration $u_{(r)}=u+\Delta_{i_1}+\cdots +\Delta_{i_r}$
is $I_r$-forbidden, where $I_r=I\setminus\{i_1,\dots,i_r\}$.
In particular, $u_{(k)}$ is $I_k$-forbidden, where $I_k=\emptyset$, 
which is impossible.  This shows that the configuration $u$ is allowed.
\end{proof}

%\begin{conjecture} {\rm \cite{Dhar}} \  
%A configuration is recurrent if and only if it 
%is stable and allowed.
%\label{conj:Dhar}
%\end{conjecture}

Dhar suggested that a configuration is recurrent if and only if it 
is stable and allowed.  
Gabrielov~\cite[Section~3, Appendix~E]{Gab1} showed that this statement 
is not true in general, and proved the conjecture for a toppling matrix 
$\Delta$ with non-negative column sums and, in particular, for 
a symmetric toppling matrix~$\Delta=L_G$
corresponding to an undirected graph $G$.
% see~\cite[Eq.~(21), App.~E]{Gab1}.
For symmetric $\Delta=L_G$, Dhar's conjecture
%~\ref{conj:Dhar} 
was also proved by Ivashkevich and Priezzhev~\cite{IP}, and 
recently by Meester, Redig, Znamenski~\cite[Theorem~5.4]{MRZ}, and
by Cori, Rossin, and Salvy~\cite[Theorem~15]{CRS}. 

The following two claims show how $G$-parking functions from
Section~\ref{sec:G-parking} are related to the sandpile model.
For a vector $u=(u_1,\dots,u_n)$, let
$u^\vee=(u_1^\vee,\dots,u_n^\vee)$, where $u_i^\vee = \Delta_{ii} - 1 -
u_i$.

\begin{lemma}
Let $G$ be a digraph with at least one oriented spanning tree, 
and let $\Delta=L_G^T$ 
be the transpose of the truncated Laplace matrix
for the digraph $G$.  
For the sandpile model associated with the toppling matrix $\Delta$,
a configuration $u$ is stable and allowed if and only if 
$u^\vee$ is a $G$-parking function.
\label{lem:Delta-G}
\end{lemma}

\begin{proof}
Parts 1 and 2 of Proposition~\ref{prop:toppling-properties} imply that
$\Delta=L_G^T$ is a toppling matrix.  
The statement of the lemma is immediate from the
definitions of allowed configurations and $G$-parking functions.
\end{proof}

A toppling matrix $\Delta$ is the transpose $L_G^T$ of 
the truncated Laplace matrix for some digraph $G$ if and only if 
it has non-negative column sums:
$$
\sum_i \Delta_{ij}\geq 0\textrm{ for any } j.
$$

Theorem~\ref{th:G-parking=trees} recovers Gabrielov's result
on recurrent configurations.

\begin{corollary} {\rm \cite[Eq.~(21)]{Gab1}} \ 
For a toppling matrix $\Delta$ with non-negative column sums, a
configuration is recurrent if and only if it is stable and allowed.

Equivalently, a configuration $u$ is recurrent if and only if $u^\vee$ is a
$G$-parking function, for $G$ and $\Delta=L_G^T$ such as in 
Lemma~\ref{lem:Delta-G}.
\label{cor:sandpiles-parking} 
\end{corollary}

\begin{proof}
Theorem~\ref{th:G-parking=trees} and the Matrix-Tree Theorem
imply that the number of
stable allowed configurations equals to $\det L_G = \det \Delta$.
According to Theorem~\ref{th:order-sandpiles},  the number of recurrent 
configurations is also equal to $\det\Delta$. 
These facts, together with Proposition~\ref{prop:recurrent->allowed},
imply the statement.
\end{proof}

Remark that we need to impose the {\it transpose\/} of Dhar's physical
conditions~(\ref{eq:toppling=laplace}) on the toppling matrix $\Delta$
in Corollary~\ref{cor:sandpiles-parking}.  The number of recurrent
configurations for a toppling matrix $\Delta$ is equal to the number of
recurrent configurations for the transposed toppling matrix $\Delta^T$,
because $\det \Delta = \det \Delta^T$.  It would be interesting to
present an explicit bijection between these two sets of configurations.


\newpage

\begin{thebibliography}{PSS2}

\bibitem[BPS]{BPS} {\sc D.~Bayer, I.~Peeva, B.~Sturmfels:}
Monomial resolutions, {\it Mathematical Research Letters \bf 5} (1998), 31--46.

\bibitem[BaSt]{BaSt} {\sc D.~Bayer, B.~Sturmfels:}
Cellular resolutions of monomial modules,
{\it Journal f\"ur die Reine und Angewandte Mathematik 
\bf 502} (1998), 123--140.

%\bibitem[CoRo]{CR} {\sc R.~Cori, D.~Rossin:}
%On the sandpile group of dual graph, {\it European Journal
%of Combinatorics \bf 21} (2000), no.~4, 447--459.

\bibitem[CRS]{CRS} {\sc R.~Cori, D.~Rossin, B.~Salvy:}
Polynomial ideals for sandpiles and their Gr\"obner bases,
{\it Theoretical Computer Science \bf 276} (2002), no.~1-2,
1--15.

\bibitem[Dhar]{Dhar} {\sc D.~Dhar:}  Self-organised critical state
of the sandpile automaton models, {\it Physical Review Letters \bf 64}
(1990), no.~14, 1613--1616.


\bibitem[Gab1]{Gab1} {\sc A.~Gabrielov:} Abelian avalanches and 
Tutte polynomials, {\it Physica A \bf 195} (1993), 253--274. 

\bibitem[Gab2]{Gab2} {\sc A.~Gabrielov:} Asymmetric abelian
avalanches and sandpile, preprint 93-65, MSI, Cornell University, 
1993.

\bibitem[IvPr]{IP} {\sc E.~V.~Ivashkevich, V.~B.~Priezzhev:}
Introduction to the sandpile model, {\it Physica A \bf 254}
(1998), 97--116.

\bibitem[Kac]{Kac}   {\sc V.~G.~Kac:} {\it Infinite dimensional Lie algebras},
third edition, Cambridge University Press, 1990.

\bibitem[Krew]{Krew} {\sc G.~Kreweras:} Une famille de polyn\^omes ayant 
plusieurs propri\'et\'es \'enumeratives, 
{\it  Periodica Mathematica Hungarica \bf 11} (1980), no.~4, 309--320.

\bibitem[MRZ]{MRZ} {\sc R.~Meester, F.~Redig, D.~Znamenski:}
The abelian sandpile: a mathematical introduction,
{\it Markov Process.\ Related Fields \bf 7} (2001), no.~4, 509--523.


\bibitem[MSY]{MSY} {\sc E.~Miller, B.~Sturmfels, K.~Yanagawa:}
Generic and cogeneric monomial ideals,
{\it Journal of Symbolic Computation \bf 29} (2000), 
no.~4-5, 691--708.

%\bibitem[Naru]{Naru} {\sc H.~Narushima:} Principle of inclusion-exclusion on
%semilattices, {\it Journal of Combinatorical Theory, Ser.~A \bf 17}
%(1974), 196--203. 

\bibitem[Naru]{Naru} {\sc H.~Narushima:} Principle of inclusion-exclusion on
partially ordered sets, {\it Discrete Mathematics \bf 42}
(1982), 243--250. 


\bibitem[PP]{PP} {\sc I.~M.~Pak, A.~E.~Postnikov:}
Resolutions for $S_n$-modules corresponding to skew hooks, 
and combinatorial applications,  
{\it  Functional Analysis and its Applications \bf 28} (1994), no.~2, 132--134.

\bibitem[PiSt]{PiSt} {\sc J.~Pitman, R.~Stanley:}  A polytope related to 
empirical distributions, plane trees, parking functions, and the associahedron,
{\it Discrete and Computational Geometry \bf 27} (2002), 603--634. 

\bibitem[PSS1]{PSS} {\sc A.~Postnikov, B.~Shapiro,  M.~Shapiro:}
Algebras of curvature forms on homogeneous manifolds,
in ``Differential Topology, Infinite-Dimensional Lie
Algebras, and Applications: D.~B.~Fuchs 60th Anniversary Collection,''
{\it AMS Translations, Ser.~2 \bf 194} (1999), 227--235.

\bibitem[PSS2]{PSS0}  {\sc A.~Postnikov, B.~Shapiro,  M.~Shapiro:}
Chern forms on flag manifolds and forests, 
{\it Proceedings of the 10-th International Conference 
on Formal Power Series and Algebraic Combinatorics,}
FPSAC'98, Fields Institute, Toronto, 1998.

\bibitem[Sche]{Schenck} {\sc H.~Schenck:}
Linear series on a special rational surface, preprint dated April~5, 2002.

\bibitem[ShSh]{SS} {\sc B.~Shapiro, M.~Shapiro:} 
On algebra generated by
Bott-Chern 2-forms on $\mathrm{SL}_n/\mathrm{B}$,
{\it C.~R.\ Acad.\ Sci.\ Paris S\'er.\ I Math.\ \bf 326}
(1998), is.~1, 75--80.

\bibitem[Sta1]{EC1}{\sc R.~P.~Stanley:}  {\it Enumerative Combinatorics,
Volume 1,} Cambridge Studies in Advanced Mathematics {\bf 49}, 
Cambridge University Press, Cambridge, 1997.

\bibitem[Sta2]{EC2}{\sc R.~P.~Stanley:}  {\it Enumerative Combinatorics,
Volume 2,} Cambridge Studies in Advanced Mathematics {\bf 62}, 
Cambridge University Press, Cambridge, 1999.

%\bibitem{St} {\sc R.~P.~Stanley:}  Hyperplane arrangements,
%parking functions and tree inversions,  Mathematical essays in honor of
%Gian-Carlo Rota, (Cambridge, MA, 1996),
%{\it Progr. Math.,\ \bf 161}  Birkh\"auser Boston, MA, (1998), 359--375.

\bibitem[Yan]{Yan} {\sc C.~H.~Yan:}  
On the enumeration of generalized parking
functions, Proceedings of the 31-st Southeastern International
Conference on Combinatorics, Graph Theory and Computing 
(Boca Raton, FL, 2000),
{\it Congressus Numerantium \bf 147} (2000), 201--209.

\end{thebibliography}
\end{document}